\documentclass[a4paper, 12pt, onecolumn]{article}
\usepackage{geometry}

\usepackage{amsmath}
\usepackage{amsthm}
\usepackage{amsfonts}
\usepackage{bbm}
\usepackage{CJK}
\usepackage{fancyhdr}
\usepackage{graphicx}
\usepackage{geometry}
\usepackage{psfrag}
\usepackage{amsfonts,amsmath,amsthm, amssymb}
\usepackage{latexsym, euscript, epic, eepic}
\usepackage{time}
\usepackage{txfonts}
\usepackage{colortbl}
\usepackage{stmaryrd}
\usepackage{mathrsfs}
\usepackage{txfonts}
\usepackage{amsfonts}
\usepackage{color}
\usepackage{lineno}
\usepackage[square, comma, sort&compress, numbers]{natbib}
\usepackage{indentfirst,latexsym,bm}

 \numberwithin{equation}{section}

\begin{document}
\newtheorem{theorem}{Theorem}[section]
\newtheorem{proposition}[theorem]{Proposition}
\newtheorem{remark}[theorem]{Remark}
\newtheorem{corollary}[theorem]{Corollary}
\newtheorem{definition}{Definition}[section]
\newtheorem{lemma}[theorem]{Lemma}
\newtheorem{conjecture}[theorem]{Conjecture}
\newtheorem{question}[theorem]{Question}
\newtheorem{problem}[theorem]{Open Problem}
\newcommand{\IN}{\mathbb{N}}
\newcommand{\IR}{\mathbb{R}}
\newcommand{\IC}{\mathbb{C}}
\newcommand{\e}{\mbox{e}}
\newcommand{\arth}{{\rm{arth}}}
\newcommand{\K}{\mathscr{K}}
\newcommand{\E}{\mathscr{E}}
\newcommand{\M}{\mathscr{M}}

\title{\large On a Conjecture Concerning the Approximates of Complete \\Elliptic Integral of the First Kind by Inverse Hyperbolic Tangent
}
\author{\normalsize Song-Liang Qiu\footnote{Corresponding author.}\,, Qi Bao, Xiao-Yan Ma, Hong-Biao Jiang
}
\date{}
\maketitle

\fontsize{12}{22}\selectfont\small

{\bf Abstract:}. Let $\K$ be the complete elliptic integral of the first kind. In this paper, the authors prove that the function $r\mapsto r^{-2}\{[\log(2\K(r)/\pi)]/\log((\arth r)/r)-3/4\}$ is strictly increasing from $(0,1)$ onto $(1/320,1/4)$, so that $[(\arth r)/r]^{3/4+r^2/320}<2\K(r)/\pi<[(\arth r)/r]^{3/4+r^2/4}$ for $r\in(0,1)$, in which all the coefficients of the exponents of the two bounds are best possible, thus proving a conjecture raised by Alzer and Qiu to be true, and giving better bounds of $\K(r)$ than those they conjectured and put in an open problem. Some other analytic properties of the complete elliptic integrals, including other kind of approximates for $\K(r)$, are obtained, too.
\\[10pt]
\\[10pt]

\section{\normalsize Introduction}\label{Sec1}

Throughout this paper, we let $\IN$ ($\arth$) denote the set of positive integers (the inverse hyperbolic tangent, respectively), $\IN_0=\IN\cup\{0\}$,  and let $r^{\thinspace\prime}=\sqrt{1-r^2}$ for each $r\in[0,1]$. For $r\in(0,1)$, {\it the complete elliptic integrals of the first and second kinds} are defined as
\begin{align}
&\begin{cases}\label{K}
\K=\K(r)=\int_0^{\pi/2}\left(1-r^2\sin^2t\right)^{-1/2}dt,\\
\K'=\K'(r)=\K(r^{\thinspace\prime}),\\
\K(0)=\pi/2,\K(1)=\infty,
\end{cases}
\end{align}
and
\begin{align}
&\begin{cases}\label{E}
\E=\E(r)=\int_0^{\pi/2}\left(1-r^2\sin^2t\right)^{1/2}dt,\\
\E'=\E'(r)=\E(r^{\thinspace\prime}),\\
\E(0)=\pi/2,\E(1)=1,
\end{cases}
\end{align}
respectively, which have wide applications and are the particular cases of the Gaussian hypergeometric function
\begin{align}\label{F}
F(a,b;c;x)={}_2F_1(a,b;c;x)=\sum_{n=0}^{\infty}\frac{(a)_n(b)_n}{(c)_nn!}x^n ~(|x|<1)
\end{align}
for complex numbers $a, b$ and $c$ with $c\neq0, -1,-2,...$, where $(a)_n=a(a+1)(a+2)\cdots(a+n-1)=\Gamma(a+n)/\Gamma(a)$ for $n\in\IN$, and $(a)_0=1$, while $\Gamma(x)=\int_0^{\infty}t^{x-1}e^{-t}dt$ is the Euler gamma function. As a matter of fact, it is well known that
\begin{align}\label{K-E}
\K(r)=\frac{\pi}{2}F\left(\frac12,\frac12;1;r^2\right) \mbox{~and~}
\E(r)=\frac{\pi}{2}F\left(-\frac12,\frac12;1;r^2\right)
\end{align}
(cf. \cite[17.3.9--17.3.10]{AS}). Moreover, by \cite[15.1.4]{AS},
\begin{align}\label{arth}
\frac{\arth r}{r}=F\left(\frac12,1;\frac32;r^2\right)=\sum_{n=0}^{\infty}\frac{1}{2n+1}r^{2n}.
\end{align}
As we known, numerous properties of the complete elliptic integrals have been obtained (cf. \cite{Ask2,AS,AAR,AQ,AQVa,AQVV,AVVb,AVV1,AVV2, Be2,B,BB,C,OLB,QMC2,QVu2,CWQ,HQM,PV,QVa,WZC} and bibliographies therein).

One of the topics of studying the properties of the complete elliptic integrals is to approximate $\K(r)$ by means of the function $r\mapsto (\arth r)/r$. Following this topic, G.D. Anderson et al. proved in \cite[Theorem 3.10]{AVV2} that the function $r\mapsto r\K(r)^2/\arth r$ ($r\mapsto r\K(r)/\arth r$) is strictly increasing (decreasing) from $(0,1)$ onto $(\pi^2/4,\infty)$ ($(1,\pi/2)$, respectively), so that
\begin{align}\label{AVVIneq1}
\max\left\{\left(\frac{\arth r}{r}\right)^{1/2}, \frac{2}{\pi}\frac{\arth r}{r}\right\}<\frac{2}{\pi}\K(r)<\frac{\arth r}{r}
\end{align}
for $r\in(0,1)$. In \cite[Theorem 18]{AQ}, Alzer and Qiu improved (\ref{AVVIneq1}) to the following double inequality
\begin{align}\label{AQIneq1}
\frac{\pi}{2}\left(\frac{\arth r}{r}\right)^{\alpha_1}<\K(r)<\frac{\pi}{2}\left(\frac{\arth r}{r}\right)^{\beta_1}
\end{align}
with best possible exponents $\alpha_1=3/4$ and $\beta_1=1$, and in \cite[Section 4]{AQ}, they posed the following
\begin{problem}\label{AQOP} What is the largest (smallest) constant $\alpha^{\ast}$ ($\beta^{\ast}$, respectively) such that the double inequality
\begin{align}\label{AQIneq2}
\frac{\pi}{2}\left(\frac{\arth r}{r}\right)^{3/4+\alpha^{\ast} r}<\K(r)<\frac{\pi}{2}\left(\frac{\arth r}{r}\right)^{3/4+\beta^{\ast} r}
\end{align}
holds for all $r\in(0,1)$?
\end{problem}

The double inequality (\ref{AQIneq2}) is obviously equivalent to
\begin{align}\label{AQIneq3}
\alpha^{\ast}<\frac{1}{r}\left[G(r)-\frac34\right]<\beta^{\ast}, \mbox{~ where ~} G(r)=\frac{\log(2\K(r)/\pi)}{\log((\arth r)/r)}
\end{align}
for $r\in(0,1)$. In \cite[Section 4]{AQ}, it was shown that the best possible constants $\alpha^{\ast}=0$ and $\beta^{\ast}\geq1/4$, and the following conjecture was put forward.
\begin{conjecture}\label{AQConj}
The function $G$ in (\ref{AQIneq3}) is strictly increasing and convex from $(0,1)$ onto $(3/4,1)$.
\end{conjecture}
If Conjecture \ref{AQConj} is true, then the best possible constants $\alpha^{\ast}$ and $\beta^{\ast}$ in (\ref{AQIneq2}) are 0 and 1/4, respectively.

In 2011, Y.-M. Chu et al. proved in \cite{CWQ1} that (\ref{AQIneq2}) holds for all $r\in(0,1)$ with the best possible constant $\beta^{\ast}=1/4$. In \cite{KS}, P\'al A. Kup\'an and R\'obert Sz\'asz studied Open Problem \ref{AQOP} and Conjecture \ref{AQConj}, and proved that for $r\in(0,1)$, the following double inequality
\begin{align}\label{KSIneq}
\frac{\pi}{2}\left(\frac{\arth r}{r}\right)^{3/4+r^4/200}<\K(r)<\frac{\pi}{2}\left(\frac{\arth r}{r}\right)^{3/4+r^2/4}
\end{align}
holds, which improves the result obtained in \cite{CWQ1}. In \cite[Section 1]{KS}, it was stated that ``It seems very difficult to prove the conjecture regarding the monotonicity of $G$". Here ``the conjecture" means Conjecture \ref{AQConj}.

Based on the results above-mentioned and our computation, the following question is natural and more reasonable than Open Problem \ref{AQOP}.
\begin{question}\label{QJM}
What are the best possible constants $\alpha$ and $\beta$ such that for $r\in(0,1)$,
\begin{align}\label{QJMIneq}
\frac{\pi}{2}\left(\frac{\arth r}{r}\right)^{3/4+\alpha r^2}<\K(r)<\frac{\pi}{2}\left(\frac{\arth r}{r}\right)^{3/4+\beta r^2}?
\end{align}
\end{question}
Let $\alpha$ and $\beta$ be the best possible constants such that (\ref{QJMIneq}) holds for $r\in(0,1)$, $G$ as in (\ref{AQIneq3}), and let
\begin{align}\label{DefOff}
f(r)=\frac{1}{r^2}\left[G(r)-\frac34\right] ~\mbox{ ~for ~} r\in(0,1).
\end{align}
Then it is easy to obtain the limiting values $f(0^+)=1/320$ and $f(1^-)=1/4$ (see Theorem \ref{Th2}(3)). It is clear that (\ref{QJMIneq}) is equivalent to
\begin{align}\label{AQIneq4}
\alpha=\inf_{0<r<1}f(r)<f(r)<\sup_{0<r<1}f(r)=\beta
\end{align}
for $r\in(0,1)$. Observe that (\ref{AQIneq2}), with the best possible constants $\alpha^{\ast}=0$ and $\beta^{\ast}=1/4$, is weaker than (\ref{QJMIneq}) if $0<\alpha<\beta\leq1/4$ in (\ref{AQIneq4}). From this point of view, the convexity of $G$ in Conjecture \ref{AQConj} is no longer so important, if our following conjecture is true.
\begin{conjecture}\label{QConj}
The function $f$ defined by (\ref{DefOff}) is strictly increasing from $(0,1)$ onto $(1/320,1/4)$.
\end{conjecture}
Clearly, the truth of Conjecture \ref{QConj} implies that the function $G$ in Conjecture \ref{AQConj}
is strictly increasing from $(0,1)$ onto $(3/4,1)$ , and gives the answer to Question \ref{QJM} (see Theorem \ref{Th2}(3)). Naturally, the proof of Conjecture \ref{QConj} is much more difficult than that of the monotonicity of $G$.

In addition, based on \cite[Theorem 3.10]{AVV2}, the following question is natural.
\begin{question}\label{Q1}
Is the function $r\mapsto r\K(r)/\arth r$ concave on $(0,1)$?
\end{question}

The main purpose of this paper is to prove that Conjecture \ref{QConj} is true, to give the complete answers to Questions \ref{QJM} and \ref{Q1} (see Theorems \ref{Th1}--\ref{Th3} stated at the end of this section), and to present some other related properties of the complete elliptic integrals of the first and second kinds, including the sharp approximates of $\K(r)$ by certain combinations in terms of $(\arth r)/r$ and polynomials (see Theorems \ref{Th4}--\ref{Th5} proved in Section \ref{Sec3}), by showing the monotonicity and convexity properties of certain combinations defined in terms of $\K(r)$, $\E(r)$ and $\arth r$. In order to prove these results, we shall establish several lemmas in Section \ref{Sec2}. We shall put forward three conjectures and an open problem in Section \ref{Sec5}, for further related studies.

\smallskip
In the sequel, we always let $a_n=[(1/2)_n/n!]^2$ and $\tilde{a}_n=3/[4(2n+1)]-a_n$ for $n\in\IN_0$, and let
\begin{align}
F_0(x)&=F\left(\frac12,\frac12;1;x\right), ~~F_1(x)=F\left(\frac12,1;\frac32;x\right),\label{F0F1}\\
G_0(x)&=F\left(\frac12,\frac12;2;x\right) \mbox{~ and ~} G_1(x)=F\left(\frac12,1;\frac52;x\right)\label{G0G1}
\end{align}
for $x\in(0,1)$. By \cite[15.2.1\&15.3.3]{AS}, we have the following relations
\begin{align}\label{F0'F1'}
F_0'(x)=\frac{G_0(x)}{4(1-x)} \mbox{~ and ~} F_1'(x)=\frac{G_1(x)}{3(1-x)}.
\end{align}

Now we state some of our main results, which will be proved in Section \ref{Sec4}.

\begin{theorem}\label{Th1}
For each $c\in(0,\infty)$, define the function $g$ on $(0,1)$ by
$$g(r)=\left(\frac34+cr^2\right)\log\frac{\arth r}{r}-\log\frac{2\K(r)}{\pi}.$$
Then $g$ is strictly increasing (decreasing) on $(0,1)$ if and only if $c\geq1/4$ ($c\leq1/320$, respectively), with
\begin{align*}
g(0^+)=0 \mbox{~ and~ } g(1^-)=
\begin{cases}
\log(\pi/2), \mbox{~if} ~c=1/4,\\
\infty, ~~~~~~~~~~~\mbox{if} ~c>1/4,\\
-\infty, ~~~~~~~~\,\mbox{if} ~c<1/4.
\end{cases}
\end{align*}
Moreover, if $c\geq1/4$ ($c\leq1/320$), then $g$ is convex (concave, respectively) on $(0,1)$.
\end{theorem}

\begin{theorem}\label{Th2}
(1) The function
$$g_1(r)\equiv\frac{1}{r^2}\left[\frac{3+r^2}{4}\log\frac{\arth r}{r}-\log\frac{2\K(r)}{\pi}\right]$$
is strictly increasing from $(0,1)$ onto $(0,\log(\pi/2))$, while the function
$$g_2(r)\equiv\frac{g_1(r)}{\log((\arth r)/r)}=\frac{1}{r^2\log((\arth r)/r)}\left[\frac{3+r^2}{4}\log\frac{\arth r}{r}-\log\frac{2\K(r)}{\pi}\right]$$
is strictly decreasing from $(0,1)$ onto $(0,79/320)$.

(2) Define the functions $g_3$ and $g_4$ on $(0,1)$ by
\begin{align*}
g_3(r)&=\frac{1}{r^2\log((\arth r)/r)}\left[\log\frac{2\K(r)}{\pi}-\left(\frac34+\frac{1}{320}r^2\right)\log\frac{\arth r}{r}\right],\\
g_4(r)&=\frac{1}{r^4}\left[\log\frac{2\K(r)}{\pi}-\left(\frac34+\frac{1}{320}r^2\right)\log\frac{\arth r}{r}\right].
\end{align*}
Then $g_3$ and $g_4$ are both strictly increasing on $(0,1)$, with ranges $(0,79/320)$ and $(0,\infty)$, respectively.

(3) The function $f$ defined by (\ref{DefOff}) is strictly increasing from $(0,1)$ onto $(1/320,1/4)$, namely, Conjecture \ref{QConj} is true. In particular, the function
$$G(r)\equiv\frac{\log(2\K(r)/\pi)}{\log((\arth r)/r)}$$
defined as in (\ref{AQIneq3}) is strictly increasing from $(0,1)$ onto $(3/4,1)$, and the double inequality
\begin{align}\label{Ineq1}
\frac{\pi}{2}\left(\frac{\arth r}{r}\right)^{3/4+\alpha r^2}<\K(r)<\frac{\pi}{2}\left(\frac{\arth r}{r}\right)^{3/4+\beta r^2}
\end{align}
holds for $r\in(0,1)$ if and only if $\alpha\leq1/320$ and $\beta\geq1/4$. Moreover, for $r\in(0,1)$,
\begin{align}\label{Ineq2}
\left(\frac{\pi}{2}\right)^{1-r^2}\left(\frac{\arth r}{r}\right)^{3/4+r^2/4}<\K(r)<\frac{\pi}{2}\left(\frac{\arth r}{r}\right)^{3/4+r^2/4}.
\end{align}
\end{theorem}

\begin{theorem}\label{Th3}
(1) The function $h_1(r)\equiv r\K(r)/\arth r$ is strictly decreasing and concave from $(0,1)$ onto $(1,\pi/2)$.

(2) For each $n\in\IN$, define the function $h_2$ on $(0,1)$ by
$$h_2(r)=\frac{1}{r^{2(n+1)}}\left[\frac{3\arth r}{4r}-\frac{2}{\pi}\K(r)-\sum_{k=0}^n\tilde{a}_k r^{2k}\right].$$
Then the coefficients of the Maclaurin series of $h_2$ are all positive. In particular, $h_2$ is strictly absolutely monotone on $(0,1)$, namely, $h_2^{(k)}(r)>0$ for all $k\in\IN_0$ and $r\in(0,1)$, with $h_2(0^+)=\tilde{a}_{n+1}>0$ and $h_2(1^-)=\infty$.

(3) The function $h_3(r)\equiv rh_2(r)/\arth r$ is strictly increasing from $(0,1)$ onto $(\tilde{a}_{n+1}, (3/4)-2/\pi)$.

(4) The function
$$h_4(r)\equiv \frac{1}{r^6}\left[\frac{\pi}{2}\left(1-\frac{1}{12}r^2-\frac{91}{2880}r^4\right)-\frac{r\K(r)}{\arth r}\right]$$
is strictly increasing from $(0,1)$ onto $(871\pi/96768,2549\pi/5760-1)$.

(5) For each $n\in\IN$, and for all $r\in(0,1)$,
\begin{align}
\frac14+\frac{3\arth r}{4r}-\sum_{k=2}^{n+1}\tilde{a}_kr^{2k}-P_1(r)r^{2n+2}
&<\frac{2}{\pi}\K(r)<\frac14+\frac{3\arth r}{4r}-\sum_{k=2}^{n+1}\tilde{a}_kr^{2k}-P_2(r)r^{2n+2},\label{K-arth1}\\
\left[P_3(r)-\alpha r^6\right]\frac{\arth r}{r}&<\frac{2}{\pi}\K(r)<\left[P_3(r)-\beta r^6\right]\frac{\arth r}{r}\label{K-arth2}
\end{align}
and
\begin{align}\label{K-arth3}
\frac14+\frac34\left(1-\delta r^4\right)\frac{\arth r}{r}&<\frac{2}{\pi}\K(r)<\frac14+\frac34\left(1-\eta r^4\right)\frac{\arth r}{r},
\end{align}
where
\begin{align*}
P_1(r)&=\frac{3\pi-8}{4\pi}\frac{\arth r}{r}-\tilde{a}_{n+1}, ~
P_2(r)=\tilde{a}_{n+1}\left(\frac{\arth r}{r}-1\right),\\
P_3(r)&=1-\frac{1}{12}r^2-\frac{91}{2880}r^4, ~\alpha=\frac{2549}{2880}-\frac{2}{\pi}=0.2484496\cdots,
\end{align*}
$\beta=871/48384=0.0180018\cdots$, $\delta=1-8/(3\pi)=0.151173\cdots$ and $\eta=1/80=0.0125$. Moreover, the constants $\alpha$, $\beta$, $\delta$ and $\eta$ in (\ref{K-arth2}) and (\ref{K-arth3}) are all best possible.
\end{theorem}

\section{\normalsize Preliminaries}\label{Sec2}

In order to prove our theorems stated in Section \ref{Sec1}, we shall prove four lemmas in this section. Firstly, let us recall the following well-known formulas: For $r\in(0,1)$,
\begin{align}
\frac{d\K}{dr}&=\frac{\E(r)-r^{\thinspace\prime2}\K(r)}{rr^{\thinspace\prime2}}, ~\frac{d\E}{dr}=\frac{\E(r)-\K(r)}{r},\label{Form1}\\
\K(r)&=\log\frac{4}{r^{\thinspace\prime}}+O\left(r^{\thinspace\prime2}\log r^{\thinspace\prime}\right) ~(r\to1) \mbox{~ and}\label{Form2}\\
&\frac{\E(r)-r^{\thinspace\prime2}\K(r)}{r^2}=\frac{\pi}{4}\sum_{n=0}^{\infty}\frac{a_n}{n+1}r^{2n}\label{Form3}
\end{align}
(see \cite[15.3.10]{AS}, \cite[Theorem 4.1]{AQVV} and \cite[(5.3)]{AQVV}, respectively), which will be frequently applied later.

Our first lemma presents several properties of the elementary function $\arth$ and the special function $\K$, most of which will be applied in the proofs of our main results.

\begin{lemma}\label{Lem1}
(1) The function $f_1(r)\equiv (r-r^{\thinspace\prime2}\arth r)/r^3$ is strictly increasing and convex from $(0,1)$ onto $(2/3,1)$.

(2) The function $f_2(r)\equiv (r-r^{\thinspace\prime2}\arth r)/(r^2r^{\thinspace\prime2}\arth r)$ is strictly increasing from $(0,1)$ onto $(2/3,\infty)$.

(3) The function $f_3(r)\equiv (r-r^{\thinspace\prime2}\arth r)/[(1+r^2)\arth r-r]$ is strictly decreasing from $(0,1)$ onto $(0,1/2)$.

(4) The function $f_4(r)\equiv[(1+r^2)\arth r-r]/[r^3\K(r)]$ is strictly increasing from $(0,1)$ onto $(8/(3\pi),2)$.

(5) The function $f_5(r)\equiv [1+(r^{\thinspace\prime2}\arth r)/r]/r^{\thinspace\prime}$ is strictly increasing and convex from $(0,1)$ onto $(2,\infty)$.

(6) The function $f_6(r)\equiv (r/r^{\thinspace\prime}-r^{\thinspace\prime}\arth r)/(\arth r-r)$ is strictly increasing from $(0,1)$ onto $(2,\infty)$.

(7) The function $f_7(r)\equiv [\log((\arth r)/r)]/r^2$ is strictly increasing from $(0,1)$ onto $(1/3,\infty)$.

(8) The functions
$$f_8(r)\equiv \frac{r^{\thinspace\prime2}\arth r}{r-r^{\thinspace\prime2}\arth r}\log\frac{\arth r}{r} \mbox{~ and~ } f_9(r)\equiv\frac{r^{\thinspace\prime2}\K(r)}{\E(r)-r^{\thinspace\prime2}\K(r)}\log\frac{2\K(r)}{\pi}$$
are both strictly decreasing from $(0,1)$ onto $(0,1/2)$.

(9) For all $r\in(0,1)$,
\begin{align}
f_{10}(r)&\equiv\frac{1}{r^4}\left[72-126r^2+19r^4-\left(72-150r^2-253r^4+167r^6\right)\frac{\arth r}{r}\right]>\frac{25688}{105},\label{f10}\\
f_{11}(r)&\equiv\frac{1}{r^2}\left[\left(13-9r^2\right)\frac{\arth r}{r}-13+6r^2\right]>\frac{14}{15}\label{f11}
\end{align}
and
\begin{align}\label{f12}
f_{12}(r)&\equiv\frac{1}{r^2}\left[\left(1-5r^2\right)\frac{\arth r}{r}-1+2r^2\right]<-\frac83.
\end{align}

(10) The functions
\begin{align*}
f_{13}(r)&\equiv\frac{1}{r^6}\left[\left(15-12r^2+r^4\right)\frac{\arth r}{r}-15+7r^2\right] \mbox{and}\\
f_{14}(r)&\equiv\frac{1}{r^4}\left[3 r^{\thinspace\prime2}-\left(3-4r^2-r^4\right)\frac{\arth r}{r}\right]
\end{align*}
are both absolutely monotone on $(0,1)$, with ranges $(8/105,\infty)$ and $(26/15,\infty)$, respectively.
\end{lemma}

{\it Proof.} (1) By (\ref{arth}), $f_1(r)$ can be written as
\begin{align}\label{f1}
f_1(r)=\frac{1-(r^{\thinspace\prime2}\arth r)/r}{r^2}=2\sum_{n=0}^{\infty}\frac{r^{2n}}{(2n+1)(2n+3)},
\end{align}
yielding the result for $f_1$.

(2) Since $f_2(r)=f_1(r)[r/(r^{\,\prime2}\arth r)]$, part (2) follows from part (1) and \cite[Lemma 3(1)]{QVV1}.

(3) Clearly, $f_3(1^-)=0$. By (\ref{arth}),
\begin{align}\label{f3}
f_3(r)&=\frac{1-(r^{\thinspace\prime2}\arth r)/r}{(1+r^2)(\arth r)/r-1}\nonumber\\
&=\frac12\left[\sum_{n=0}^{\infty}\frac{r^{2n}}{(2n+1)(2n+3)}\right]\left[\sum_{n=0}^{\infty}\frac{(n+1)r^{2n}}{(2n+1)(2n+3)}\right]^{-1}.
\end{align}
Applying \cite[Lemma 2.1]{PV}, one can easily obtain the monotonicity of $f_3$. By (\ref{f3}), $f_3(0^+)=1/2$.

(4) By (\ref{K-E}) and (\ref{arth}), we can write $f_4(r)$ as
\begin{align}\label{f4}
f_4(r)&=\frac{r^{-3}[(1+r^2)\arth r-r]}{\K(r)}\nonumber\\
&=\frac{8}{\pi}\left[\sum_{n=0}^{\infty}\frac{n+1}{(2n+1)(2n+3)}r^{2n}\right]\left(\sum_{n=0}^{\infty}a_nr^{2n}\right)^{-1}.
\end{align}
Let $\lambda_n=(n+1)/[(2n+1)(2n+3)a_n]$ for $n\in\IN_0$. Then
$$\frac{\lambda_{n+1}}{\lambda_n}=\frac{(n+2)(2n+1)a_n}{(n+1)(2n+5)a_{n+1}}=\frac{4n^2+12n+8}{4n^2+12n+5}>1,$$
which shows that the sequence $\{\lambda_n\}$ is strictly increasing, and hence by \cite[Lemma 2.1]{PV}, $f_4$ is strictly increasing on $(0,1)$. By (\ref{f4}), $f_4(0^+)=8/(3\pi)$. Since
\begin{align}\label{lim}
\lim_{r\to1}r\K(r)/\arth r=1
\end{align}
by \cite[Theorem 3.10]{AVV2} or (\ref{Form2}), we obtain the limiting value $f_4(1^-)=2$.

(5) The limiting values of $f_5$ are clear. By differentiation, $f_5'(r)=rf_1(r)/r^{\thinspace\prime3}$,
and hence the monotonicity and convexity properties of $f_5$ follow from part (1).

(6) Let $G_2(r)=(r/r^{\thinspace\prime})-r^{\thinspace\prime}\arth r$ and $G_3(r)=\arth r-r$. Then $G_2(0)=G_3(0)=0$, $G_2(1^-)=G_3(1^-)=\infty$, $f_6(r)=G_2(r)/G_3(r)$, and
\begin{align}\label{G2'G3'}
G_2'(r)/G_3'(r)=f_5(r).
\end{align}
Hence the monotonicity of $f_6$ follows from \cite[Theorem 1.25]{AVVb} and part (5).

By l'H\^opital's rule and (\ref{G2'G3'}), $f_6(0^+)=2$ and $f_6(1^-)=\infty$.

(7) Let $G_4(r)=\log((\arth r)/r)$ and $G_5(r)=r^2$. Then $G_4(0^+)=G_5(0)=0$, $f_7(r)=G_4(r)/G_5(r)$ and
\begin{align*}
G_4'(r)/G_5'(r)=f_2(r)/2.
\end{align*}
Hence the monotonicity of $f_7$ follows from \cite[Theorem 1.25]{AVVb} and part (2). By l'H\^opital's rule, $f_7(0^+)=f_2(0^+)/2=1/3$. Clearly, $f_7(1^-)=\infty$.

(8) Let $G_6(r)=r/(r^{\thinspace\prime2}\arth r)-1$. Then $f_8(r)=G_4(r)/G_6(r)$, $G_4(0^+)=G_6(0^+)=0$, $G_4(1^-)=G_6(1^-)=\infty$, and by differentiation,
\begin{align*}
\frac{G_4'(r)}{G_6'(r)}=f_3(r)\frac{r^{\thinspace\prime2}\arth r}{r}
\end{align*}
which is strictly decreasing from $(0,1)$ onto $(0,1/2)$ by part (3) and \cite[Lemma 3(1)]{QVV1}. Hence the result for $f_8$ follows from \cite[Theorem 1.25]{AVVb}, part (3) and l'H\^opital's rule.

Similarly, since
$$\lim_{r\to0}\log\frac{2\K(r)}{\pi}=\lim_{r\to0}\frac{\E(r)-r^{\thinspace\prime2}\K(r)}{r^{\thinspace\prime2}\K(r)}=0,$$
and since the ratio
\begin{align*}
&\left[\frac{d}{dr}\left(\log\frac{2\K(r)}{\pi}\right)\right]
\left[\frac{d}{dr}\left(\frac{\E(r)-r^{\thinspace\prime2}\K(r)}{r^{\thinspace\prime2}\K(r)}\right)\right]^{-1}\\
&=\frac{r^{\thinspace\prime2}\K(r)\left[\E(r)-r^{\thinspace\prime2}\K(r)\right]}{\E(r)[\K(r)-\E(r)]+\K(r)\left[\E(r)-r^{\thinspace\prime2}\K(r)\right]}\\
&=\left[\frac{1}{r^{\thinspace\prime2}}+\frac{\E(r)}{r^{\thinspace\prime}}\cdot\frac{1}{r^{\thinspace\prime}\K(r)}\cdot\frac{\K(r)-\E(r)}{\E(r)-r^{\thinspace\prime2}\K(r)}\right]^{-1}
\end{align*}
is strictly decreasing on $(0,1)$ by \cite[Theorem 3.21(6)\&(8)]{AVVb}, the result for $f_9$ follows from \cite[Theorem 1.25]{AVVb} and l'H\^opital's rule.

(9) It follows from (\ref{arth}) that
\begin{align*}
f_{10}(r)=&\frac{1}{r^4}\left[19r^4-126r^2-72\sum_{n=1}^{\infty}\frac{r^{2n}}{2n+1}
+150r^2\sum_{n=0}^{\infty}\frac{r^{2n}}{2n+1}\right.\\
&\left.+253r^4\sum_{n=0}^{\infty}\frac{r^{2n}}{2n+1}-167r^6\sum_{n=0}^{\infty}\frac{r^{2n}}{2n+1}\right]\\
=&\frac{1}{r^4}\left[19r^4-72\sum_{n=2}^{\infty}\frac{r^{2n}}{2n+1}+150r^2\sum_{n=1}^{\infty}\frac{r^{2n}}{2n+1}\right.\\
&\left.+253r^4\sum_{n=0}^{\infty}\frac{r^{2n}}{2n+1}-167r^6\sum_{n=0}^{\infty}\frac{r^{2n}}{2n+1}\right]\\
=&19-\sum_{n=0}^{\infty}\frac{72 r^{2n}}{2n+5}+\sum_{n=0}^{\infty}\left(\frac{150}{2n+3}
+\frac{253}{2n+1}\right)r^{2n}-r^2\sum_{n=0}^{\infty}\frac{167r^{2n}}{2n+1}\\
=&\frac{1538}{5}-\frac{1322}{21}r^2+2r^4\sum_{n=0}^{\infty}\frac{656n^3+5540n^2+12140n+3951}{(2n+9)(2n+7)(2n+5)(2n+3)}r^{2n}\\
>&\frac{1538}{5}-\frac{1322}{21}+2r^4\sum_{n=0}^{\infty}\frac{656n^3+5540n^2+12140n+3951}{(2n+9)(2n+7)(2n+5)(2n+3)}r^{2n}\\
>&\frac{1538}{5}-\frac{1322}{21}=\frac{25688}{105},\\
f_{11}(r)=&13\sum_{n=0}^{\infty}\frac{r^{2n}}{2n+3}-9\sum_{n=0}^{\infty}\frac{r^{2n}}{2n+1}+6\\
=&\frac43-\frac25r^2+2\sum_{n=2}^{\infty}\frac{4n-7}{(2n+1)(2n+3)}r^{2n}>\frac43-\frac25r^2>\frac{14}{15}
\end{align*}
and
\begin{align*}
f_{12}(r)=&2-2\sum_{n=0}^{\infty}\frac{4n+7}{(2n+1)(2n+3)}r^{2n}<2-\frac{14}{3}=-\frac83.
\end{align*}

(10) Clearly, $f_{13}(1^-)=f_{14}(1^-)=\infty$. Applying (\ref{arth}), we obtain the following series expansions
\begin{align*}
f_{13}(r)=8\sum_{n=0}^{\infty}\frac{(n+1)(2n+1)}{(2n+3)(2n+5)(2n+7)}r^{2n}
\end{align*}
and
\begin{align*}
f_{14}(r)=2\sum_{n=0}^{\infty}\frac{4n^2+20n+13}{(2n+1)(2n+3)(2n+5)}r^{2n},
\end{align*}
from which the absolute monotonicity and the limiting values of $f_{13}$ and $f_{14}$ as $r\to0$ follow.
 $\Box$

\smallskip
Our next lemma presents some relations of $F_0(r^2)=2\K(r)/\pi$ and $F_1(r^2)=(\arth r)/r$.

\begin{lemma}\label{Lem2}
(1) The function $f_{15}(r)\equiv F_1(r)/F_0(r)$ is strictly increasing from $(0,1)$ onto $(1,\pi/2)$.

(2) The function $f_{16}(r)\equiv G_1(r)/G_0(r)=3 F_1'(r)/[4F_0'(r)]$ is strictly increasing from $(0,1)$ onto $(1,3\pi/8)$.

(3) The function
$$f_{17}(r)\equiv F_0(r)-F_1(r)F_0'(r)/F_1'(r)$$
is strictly increasing from $(0,1)$ onto $(1/4,(\log4)/\pi)$.
\end{lemma}

{\it Proof.} Part (1) is implied by \cite[Theorem 3.10]{AVV2}, and part (2) follows from \cite[Lemmas 3.2(2)]{QMC2}.

By (\ref{F0'F1'}), $f_{17}(r)=F_0(r)-3F_1(r)G_0(r)/[4G_1(r)]$. Since
\begin{align*}
f_{17}^{\,\prime}(r)=-F_1(r)\frac{d}{dr}\left[\frac{F_0'(r)}{F_1'(r)}\right]=-\frac34F_1(r)\frac{d}{dr}\left[\frac{1}{f_{16}(r)}\right]>0
\end{align*}
by part (2), the monotonicity of $f_{17}$ follows.

Clearly, $f_{17}(0)=1/4$. By (\ref{F0'F1'}),
\begin{align*}
f_{17}(r)&=-\frac{F_0'(r)}{F_1'(r)}\left[F_1(r)-F_0(r)\frac{F_1'(r)}{F_0'(r)}\right]\\
&=-\frac{3}{4f_{16}(r)}\left[F_1(r)-F_0(r)\frac{F_1'(r)}{F_0'(r)}\right].
\end{align*}
Hence the limiting value $f_{17}(1^-)=(\log4)/\pi$ follows from part (2) and \cite[Lemma 3.7(2)]{QMC2}. $\Box$

\begin{lemma}\label{Lem3}
For $r\in(0,1)$, the coefficients of the Maclaurin series of the function
$$f_{18}(r)\equiv \frac{1}{r^8}\left[\left(1-\frac12r^2-\frac{1}{16}r^4-\frac{1}{32}r^6\right)\K(r)-\E(r)\right]$$
are all \emph{}positive, with $f_{18}(0^+)=41\pi/4096$ and $f_{18}(1^-)=\infty$.
\end{lemma}

{\it Proof.} By (\ref{F}) and (\ref{K-E}), we have
\begin{align}\label{f18}
f_{18}(r)=&\frac{\pi}{2r^8}\left[\left(1-\frac12r^2-\frac{1}{16}r^4-\frac{1}{32}r^6\right)\sum_{n=0}^{\infty}a_nr^{2n}
+\sum_{n=0}^{\infty}\frac{a_n}{2n-1}r^{2n}\right]\nonumber\\
=&\frac{\pi}{2r^8}\left[\sum_{n=1}^{\infty}\left(\frac{2n}{2n-1}a_n-\frac12 a_{n-1}\right)r^{2n}-\frac{1}{16}\sum_{n=2}^{\infty}a_{n-2}r^{2n}
-\frac{1}{32}\sum_{n=3}^{\infty}a_{n-3}r^{2n}\right]\nonumber\\
=&\frac{\pi}{2r^8}\left[\sum_{n=2}^{\infty}\left(\frac{2n}{2n-1}a_n-\frac12 a_{n-1}\right)r^{2n}-\frac{1}{16}r^4
-\sum_{n=3}^{\infty}\left(\frac{1}{16}a_{n-2}+\frac{1}{32}a_{n-3}\right)r^{2n}\right]\nonumber\\
=&\frac{\pi}{2r^8}\sum_{n=4}^{\infty}\left(\frac{2n}{2n-1}a_n-\frac12 a_{n-1}-\frac{1}{16}a_{n-2}-\frac{1}{32}a_{n-3}\right)r^{2n}\nonumber\\
=&\frac{\pi}{2}\sum_{n=0}^{\infty}\left(\frac{2n+8}{2n+7}a_{n+4}-\frac12 a_{n+3}-\frac{1}{16}a_{n+2}-\frac{1}{32}a_{n+1}\right)r^{2n}\nonumber\\
=&\frac{\pi}{2}\sum_{n=0}^{\infty}b_n r^{2n},
\end{align}
where
$$b_n=\frac{2n+8}{2n+7}a_{n+4}-\frac12 a_{n+3}-\frac{1}{16}a_{n+2}-\frac{1}{32}a_{n+1}.$$
Since $4a_{n+1}=[(2n+1)/(n+1)]^2a_n$, we have
\begin{align}\label{bn}
b_n=&\left[\frac{n+4}{2(2n+7)}\left(\frac{2n+7}{n+4}\right)^2-\frac12\right]a_{n+3}-\frac{1}{32}\left(2+\frac{a_{n+1}}{a_{n+2}}\right)a_{n+2}\nonumber\\
=&\frac{n+3}{2(n+4)}a_{n+3}-\frac{6n^2+20n+17}{16(2n+3)^2}a_{n+2}\nonumber\\
=&\left[\frac{n+3}{8(n+4)}\left(\frac{2n+5}{n+3}\right)^2-\frac{6n^2+20n+17}{16(2n+3)^2}\right]a_{n+2}\nonumber\\
=&\frac{26n^4+194n^3+523n^2+601n+246}{16(n+3)(n+4)(2n+3)^2}a_{n+2}\nonumber\\
=&\frac{(n+2)(n+1)(26n^2+116n+123)}{16(n+3)(n+4)(2n+3)^2}a_{n+2}>0
\end{align}
for all $n\in\IN_0$. Hence it follows from (\ref{f18}) and (\ref{bn}) that all the coefficients of the Maclaurin series of $f_{18}$ are positive.
By (\ref{f18}), $f_{18}(0^+)=\pi b_0/2=41\pi/4096$. Clearly, $f_{18}(1^-)=\infty$. $\Box$

\bigskip
The following lemma plays an important role in the proof of Theorem \ref{Th3}.

\begin{lemma}\label{Lem4}
For $n\in\IN_0$, let $c_n=(2n+3)a_{n+1}$,
\begin{align*}
\tilde{c}_n&=\frac{(22n+83)(2n+9)(2n+7)(2n+5)(2n+3)}{(n+4)^2\left(32276n^2+123808n+110907\right)}a_{n+3},\\
d_n&=(2n+1)\left[\frac{10196n^2+39096n+34975}{(2n+7)(2n+5)(2n+3)}-2880a_{n+3}\right].
\end{align*}

(1) The sequence $\{c_n\}$ is strictly decreasing with $c_0=3/4$ and $c_{\infty}=\lim_{n\to\infty}c_n=2/\pi$.

(2) The sequence $\{\tilde{c}_n\}$ is strictly decreasing with $\tilde{c}_0=217875/50475008=0.0043164\cdots$ and $\tilde{c}_{\infty}\equiv\lim_{n\to\infty}\tilde{c}_n=88/(8069\pi)=0.0034714\cdots$.

(3) The sequence $\{d_n\}$ is strictly increasing with $d_0=4355/84=51.845238\cdots$ and $d_{\infty}\equiv\lim_{n\to\infty}d_n=2549-5760/\pi=715.535055\cdots$.

\end{lemma}

{\it Proof.} (1) Clearly, $c_0=3/4$. By applying \cite[6.1.8 \& 6.1.47]{AS}, we obtain the value $c_{\infty}=\lim_{n\to\infty}c_n=2/\pi$.

Computation gives
\begin{align*}
\frac{c_{n+1}}{c_n}&=\frac{(2n+5)a_{n+2}}{(2n+3)a_{n+1}}=\frac{2n+5}{2n+3}\left(\frac{n+3/2}{n+2}\right)^2=\frac{4n^2+16n+15}{4n^2+16n+16}<1,
\end{align*}
yielding the monotonicity of the sequence $\{c_n\}$.

(2) It is easy to see that $\tilde{c}_0=217875/50475008$. By \cite[6.1.8 \& 6.1.47]{AS}, $\lim_{n\to\infty}\tilde{c}_n=88/(8069\pi)$.

By computation, we obtain
\begin{align*}
\frac{\tilde{c}_{n+1}}{\tilde{c}_n}=&\frac{(n+4)^2(22n+105)(2n+11)\left(32276n^2+123808n+110907\right)}
{(n+5)^2(22n+83)(2n+3)\left(32276n^2+188360n+266991\right)}\frac{a_{n+4}}{a_{n+3}}\\
=&\frac{(2n+7)^2(22n+105)(2n+11)\left(32276n^2+123808n+110907\right)}
{4(n+5)^2(22n+83)(2n+3)\left(32276n^2+188360n+266991\right)}=\frac{\Delta_1(n)}{\Delta_2(n)}
\end{align*}
for $n\in\IN_0$, where
\begin{align*}
\Delta_1(n)=&5680576n^6+119909248n^5+1023083184n^4+4501649024n^3\\
&+10727855348n^2+13050014376n+6276781665 \mbox{~ and}\\
\Delta_2(n)=&5680576n^6+119909248n^5+1026984672n^4+4552315328n^3\\
&+10968049916n^2+13543585560n+6648075900.
\end{align*}
Since
\begin{align*}
\Delta_2(n)-\Delta_1(n)=&3901488n^4+50666304n^3+240194568n^2\\
&+493571184n+371294235>371294235
\end{align*}
for $n\in\IN_0$, $\tilde{c}_{n+1}<\tilde{c}_n$ for $n\in\IN_0$, namely, the sequence $\{\tilde{c}_n\}$ is strictly decreasing.

(3) Put $\Delta=d_{n+1}-d_n$. Then by computation, we obtain
\begin{align*}
\Delta=&(2n+3)\frac{10196n^2+59488n+84267}{(2n+9)(2n+7)(2n+5)}-(2n+1)\frac{10196n^2+39096n+34975}{(2n+7)(2n+5)(2n+3)}\\
&-2880\left[(2n+3)\frac{a_{n+4}}{a_{n+3}}-(2n+1)\right]a_{n+3}\\
=&4\left[\frac{32276n^2+123808n+110907}{(2n+9)(2n+7)(2n+5)(2n+3)}-180\frac{22n+83}{(n+4)^2}a_{n+3}\right],
\end{align*}
and hence
\begin{align*}
d_{n+1}>d_n&\Longleftrightarrow\frac{32276n^2+123808n+110907}{(2n+9)(2n+7)(2n+5)(2n+3)}>180\frac{22n+83}{(n+4)^2}a_{n+3}\\
&\Longleftrightarrow\frac{(22n+83)(2n+9)(2n+7)(2n+5)(2n+3)}{(n+4)^2\left(32276n^2+123808n+110907\right)}a_{n+3}<\frac{1}{180}\\
&\Longleftrightarrow \tilde{c}_n<\frac{1}{180}=0.0055555\cdots
\end{align*}
which is true by part (2). Hence the monotonicity of the sequence $\{d_n\}$ follows. Clearly, $d_0=34975/105-2880a_3=4355/84$. By \cite[6.1.8 \& 6.1.47]{AS}, we obtain the value $d_{\infty}=2549-5760/\pi$. $\Box$

%

\section{\normalsize Some Properties of the Complete Elliptic Integrals}\label{Sec3}

In this section, we present several properties of the complete elliptic integrals $\K(r)$ and $\E(r)$, which are also needed in the proofs of Theorems \ref{Th1} and \ref{Th2}.

\begin{theorem}\label{Th4}
Let $f_{18}$ be as in Lemma \ref{Lem3}. Then the function $f_{19}(r)\equiv f_{18}(r)/\K(r)$
is strictly increasing from $(0,1)$ onto $(41/2048,13/32)$. In particular, for $r\in(0,1)$,
\begin{align}\label{Ineq3}
P(r)<\E(r)/\K(r)<Q(r),
\end{align}
where $Q(r)=1-r^2/2-r^4/16-r^6/32-41 r^8/2048$ and
\begin{align*}
P(r)=1-\frac12r^2-\frac{1}{16}r^4-\frac{1}{32}r^6-\frac{13}{32}r^8=r^{\thinspace\prime2}\left(1+r^2\frac{16+14r^2+13r^4}{32}\right).
\end{align*}
\end{theorem}

{\it Proof.} For $n\in\IN_0$, let $b_n$ be as in the proof of Lemma \ref{Lem3}, and $C_n=b_n/a_n$. Then by (\ref{f18}), we have
\begin{align}\label{f19}
f_{19}(r)=\left(\sum_{n=0}^{\infty}b_nr^{2n}\right)\left(\sum_{n=0}^{\infty}a_nr^{2n}\right)^{-1}.
\end{align}
It follows from (\ref{bn}) that
\begin{align*}
C_n=&\frac{(n+2)(n+1)\left(26n^2+116n+123\right)}{256(n+3)(n+4)(2n+3)^2}\left(\frac{2n+3}{n+2}\right)^2\left(\frac{2n+1}{n+1}\right)^2\\
=&\frac{(2n+1)^2\left(26n^2+116n+123\right)}{256(n+4)(n+3)(n+2)(n+1)}
\end{align*}
and
\begin{align*}
\frac{C_{n+1}}{C_n}=&\frac{(n+1)(2n+3)^2\left(26n^2+168n+265\right)}
{(n+5)(2n+1)^2\left(26n^2+116n+123\right)}\\
=&\frac{104n^5+1088n^4+4294n^3+8002n^2+7077n+2385}{104n^5+1088n^4+3822n^3+5518n^2+3163n+615}>1
\end{align*}
which shows that the sequence $\{C_n\}$ is strictly increasing. Consequently, the monotonicity of $f_{19}$ follows from \cite[Lemma 2.1]{PV} and (\ref{f19}).

Clearly, $f_{19}(0^+)=C_0=41/2048$ and $f_{19}(1^-)=13/32$. The double inequality (\ref{Ineq3}) is clear. $\Box$

\begin{theorem}\label{Th5}
(1) Define the functions $f_{20}$ and $f_{21}$ on $(0,1)$ by
\begin{align*}
f_{20}(r)&=\left(\frac34+\frac14r^2\right)\left(r-r^{\thinspace\prime2}\arth r\right)\K(r)-\left[\E(r)-r^{\thinspace\prime2}\K(r)\right]\arth r,\\
f_{21}(r)&=\left[\E(r)-r^{\thinspace\prime2}\K(r)\right]\arth r-\left(\frac34+\frac{1}{160}r^2\right)\left(r-r^{\thinspace\prime2}\arth r\right)\K(r).
\end{align*}
Then $f_{20}$ and $f_{21}$ are both strictly increasing and convex on $(0,1)$, with ranges $(0,\log2)$ and $(0,\infty)$, respectively.

(2) The function
$$f_{22}(r)\equiv\frac{(r-r^{\thinspace\prime2}\arth r)\E(r)-2r^{\thinspace\prime2}(\arth r-r)\K(r)}{r^5}$$
is strictly increasing from $(0,1)$ onto $(\pi/30,1)$.

(3) The function
$$f_{23}(r)\equiv 2\frac{r^{\thinspace\prime2}\arth r}{r}\K(r)+\frac{r-r^{\thinspace\prime2}\arth r}{r^3}\E(r)$$
is strictly decreasing from $(0,1)$ onto $(1,4\pi/3)$.

(4) The function
$$f_{24}(r)\equiv\frac{1}{r^2}\left\{4\frac{[\E(r)-r^{\thinspace\prime2}\K(r)]\arth r}{\left(r-r^{\thinspace\prime2}\arth r\right)\K(r)}-3\right\}$$
is strictly increasing from $(0,1)$ onto $(1/40,1)$. In particular, for $r\in(0,1)$,
\begin{align}\label{f24}
\frac34+\frac{r^2}{4}P_4(r)<\frac{[\E(r)-r^{\thinspace\prime2}\K(r)]\arth r}{\left(r-r^{\thinspace\prime2}\arth r\right)\K(r)}<\frac34+\frac{r^2}{4},
\end{align}
where
$$P_4(r)=\max\left\{\frac{1}{40}, ~1-\frac{\log 16}{r\left(r-r^{\thinspace\prime2}\arth r\right)\K(r)}\right\}.$$
\end{theorem}

{\it Proof.} (1) By (\ref{Form1}) and by differentiation, we obtain
\begin{align}\label{f19'}
f_{20}^{\,\prime}(r)=&\frac12\left[2r\left(1+r^2\right)\arth r+r^2\right]\K+\frac14\left(3+r^2\right)\left(r-r^{\thinspace\prime2}\arth r\right)\frac{\E-r^{\thinspace\prime2}\K}{rr^{\thinspace\prime2}}\nonumber\\
&-r\K\arth r-\frac{\E-r^{\thinspace\prime2}\K}{r^{\thinspace\prime2}}\nonumber\\
=&\frac14\left\{\left[\left(3r^4-2r^2+3\right)\frac{\arth r}{r}+r^2+1\right]\K-\left[\left(3+r^2\right)\frac{\arth r }{r}+1\right]\E\right\}\nonumber\\
=&\frac14\left\{\left[1+\left(3+r^2\right)\frac{\arth r}{r}\right](\K-\E)+r^2\left(1-3\frac{r^{\thinspace\prime2}\arth r }{r}\right)\K\right\}\nonumber\\
=&\frac14r^2\K\left\{\left[1+\left(3+r^2\right)\frac{\arth r}{r}\right]\frac{\K-\E}{r^2\K}+1-3\frac{r^{\thinspace\prime2}\arth r }{r}\right\}\nonumber\\
=&\frac14G_7(r)G_8(r),
\end{align}
where $G_7(r)=r^2\K$ and
$$G_8(r)=\left[1+\left(3+r^2\right)\frac{\arth r}{r}\right]\frac{\K-\E}{r^2\K}+1-3\frac{r^{\thinspace\prime2}\arth r}{r}.$$

Clearly, $G_7$ is strictly increasing from $(0,1)$ onto $(0,\infty)$. By \cite[Lemma 5.2(3)]{AQVV} and \cite[Lemma 3(1)]{QVV1}, we see that $G_8$ is strictly increasing from $(0,1)$ onto $(0,\infty)$. Hence it follows from (\ref{f19'}) that $f_{20}^{\,\prime}$ is strictly increasing from $(0,1)$ onto $(0,\infty)$. This yields the monotonicity and convexity of $f_{20}$.

Clearly, $f_{20}(0^+)=0$. By l'H\^opital's rule,
\begin{align*}
\lim_{r\to1}\frac{4\E-r(3+r^2)}{r^{\thinspace\prime}}=\lim_{r\to1}\frac{1}{r}\left[4 r^{\thinspace\prime}\frac{\K-\E}{r}+3 r^{\thinspace\prime}\left(1+r^2\right)\right]=0,
\end{align*}
and hence by (\ref{Form2}), we obtain
\begin{align*}
f_{20}(1^-)=&\lim_{r\to1}\left[\frac{3+r^2}{4}\left(r-r^{\thinspace\prime2}\arth r\right)\K -\left(\E-r^{\thinspace\prime2}\K\right)\arth r\right]\\
=&\lim_{r\to1}\left[\frac{r(3+r^2)}{4}\log\frac{4}{r^{\thinspace\prime}}-\frac{\E}{2}\log\frac{1+r}{1-r}\right]\\
=&\log2+\frac18\lim_{r\to1}\frac{4\E-r(3+r^2)}{r^{\thinspace\prime}}\cdot r^{\thinspace\prime}\log(1-r)=\log2.
\end{align*}

Next, for $r\in(0,1)$, let
\begin{align*}
G_9(r)=&\left[40r-r^3+\left(120+r^2\right)r^{\thinspace\prime2}\arth r\right]\E\\
&-r^{\thinspace\prime2}\left[40r+r^3+\left(120-41r^2+3r^4\right)\arth r\right]\K,\\
G_{10}(r)=&4\left[40-\left(2r^3+79r\right)\arth r-r^2\right]\E\\
&-\left[160-\left(15r^5-131r^3+280r\right)\arth r-r^4-120r^2\right]\K,\\
G_{11}(r)=&2G_{12}(r)\K+\left[\left(47r^4+469r^2-352\right)\arth r+5r^3-368r\right]\E \mbox{~and}\\
G_{12}(r)=&r^{\thinspace\prime2}\left(30r^4-127r^2+158\right)\arth r
+6r^5-126r^3+202r.
\end{align*}
Then by differentiation,
\begin{align}
f_{21}^{\,\prime}(r)=&r\K\arth r+\frac{\E-r^{\thinspace\prime2}\K}{r^{\thinspace\prime2}}-\frac{r}{80}\left[\left(2r^2+119\right)\arth r+r\right]\K\nonumber\\
&-\frac{120+r^2}{160}\left(r-r^{\thinspace\prime2}\arth r\right)\frac{\E-r^{\thinspace\prime2}\K}{rr^{\thinspace\prime2}}
=\frac{G_9(r)}{160 rr^{\thinspace\prime2}},\label{f20'}\\
G_9'(r)=&\left[\left(18r^5-176r^3+322r\right)\arth r+2r^4+158r^2-160\right]\K\nonumber\\
&+\left[\left(3r^6-44r^4+161r^2-120\right)\arth r+r^5+39r^3-40r\right]\frac{\E-r^{\thinspace\prime2}\K}{rr^{\thinspace\prime2}}\nonumber\\
&-\left[\left(4r^3+238r\right)\arth r+2r^2-160\right]\E\nonumber\\
&+\left[\left(r^4+119r^2-120\right)\arth r+r^3-40r\right]\frac{\K-\E}{r}=G_{10}(r),\label{H8'}\\
G_{10}'(r)=&\frac{1}{r^{\thinspace\prime2}}\left[520r-367r^3+11r^5-\left(75r^6-468r^4+673r^2-280\right)\arth r\right]\K\nonumber\\
&+\left[\left(15r^5-131r^3+280r\right)\arth r+r^4+120r^2-160\right]\frac{\E-r^{\thinspace\prime2}\K}{rr^{\thinspace\prime2}}\nonumber\\
&+\frac{4}{r^{\thinspace\prime2}}\left[\left(6r^4+73r^2-79\right)\arth r-81r\right]\E\nonumber\\
&+4\left[\left(2r^3+79r\right)\arth r+r^2-40\right]\frac{\K-\E}{r}
=\frac{G_{11}(r)}{r^{\thinspace\prime2}}.\label{G9'}
\end{align}

It is clear that for all $r\in(0,1)$,
$$G_{12}(r)=r^{\thinspace\prime2}\left(158r^{\thinspace\prime2}+31r^2+30r^4\right)\arth r+r\left(202r^{\thinspace\prime2}+76r^2+6r^4\right)>0.$$
Let $f_{10}$ be as in Lemma \ref{Lem1}(9). Then by (\ref{f10}) and (\ref{Ineq3}), we have
\begin{align}\label{SignG11}
\frac{G_{11}(r)}{\E}>&\frac{4G_{12}(r)}{2-r^2}+\left(47r^4+469r^2-352\right)\arth r\nonumber\\
&+5r^3-368r=\frac{r^5}{2-r^2}f_{10}(r)>0.
\end{align}

It follows from (\ref{f10}) and (\ref{G9'})--(\ref{SignG11}) that
$$G_{10}(r)=G_9'(r)\left[\frac{d}{dr}(r)\right]^{-1}$$
is strictly increasing on $(0,1)$ with $G_{10}(0^+)=0$, and hence by \cite[Theorem 1.25]{AVVb}, the function $r\mapsto G_9(r)/r$ is strictly increasing on $(0,1)$ with $G_9(0^+)=0$. Consequently, by (\ref{f20'}), $f_{21}'$ is positive and strictly increasing on $(0,1)$. This yields the monotonicity and convexity of $f_{21}$.

Clearly, $f_{21}(0^+)=0$. By (\ref{lim}),
$$\lim_{r\to1}\left[\E(r)-r\left(\frac34+\frac{1}{160}r^2\right)\frac{\K(r)}{\arth r}\right]=\frac{39}{160},$$
so that we obtain the limiting value
\begin{align*}
f_{21}(1^-)=&\lim_{r\to1}\left[\E(r)-r\left(\frac34+\frac{1}{160}r^2\right)\frac{\K(r)}{\arth r}\right]\arth r=\infty.
\end{align*}

(2) Let $f_{11}$ and $f_{13}$ be as in Lemma \ref{Lem1}. Then by differentiation,
\begin{align}\label{f21'}
r^6f_{22}^{\,\prime}(r)=&r\left[2r\E\arth r+\left(r-r^{\,\prime2}\arth r\right)\frac{\E-\K}{r}-2r(3r-2\arth r)\K\right.\nonumber\\
&\left.-2(\arth r-r)\frac{\E-r^{\,\prime2}\K}{r}\right]-5\left[\left(r-r^{\,\prime2}\arth r\right)\E-2 r^{\,\prime2}(\arth r-r)\K\right]\nonumber\\
=&\left[\left(13-9r^2\right)\arth r-13r+6r^3\right]\K-2\left(r-r^{\,\prime2}\arth r\right)\E\nonumber\\
=&r^3f_{11}(r)\K-2\left(r-r^{\,\prime2}\arth r\right)\E.
\end{align}
By (\ref{Ineq3}) and Lemma \ref{Lem1}(10), and by (\ref{f21'}),
\begin{align*}
r^6f_{22}^{\,\prime}(r)>&r^3f_{11}(r)\frac{2\E}{2-r^2}-2\left(r-r^{\,\prime2}\arth r\right)\E\\
=&\frac{2\E}{2-r^2}\left[r^3f_{11}(r)-\left(2-r^2\right)\left(r-r^{\,\prime2}\arth r\right)\right]\\
=&\frac{2\E}{2-r^2}\left[\left(15-12r^2+r^4\right)\arth r-15r+7r^3\right]\\
=&\frac{2r^7\E}{2-r^2}f_{13}(r)>\frac{16r^7\E}{105(2-r^2)}>0
\end{align*}
for all $r\in(0,1)$, and hence the monotonicity of $f_{22}$ follows.

Clearly, $f_{22}(1^-)=1$. By Lemma \ref{Lem1}(1), (\ref{arth}) and (\ref{Form3}), we obtain
\begin{align*}
f_{22}(0^+)=&\lim_{r\to0}\frac{1}{r^5}\left\{\left(r-r^{\,\prime2}\arth r\right)\left(\E-r^{\,\prime2}\K\right)-r^{\,\prime2}\left[\left(3-r^2\right)\arth r-3r\right]\K\right\}\\
=&\lim_{r\to0}\left[\frac{(r-r^{\,\prime2}\arth r)}{r^3}\frac{\E-r^{\,\prime2}\K}{r^2}-4r^{\,\prime2}\K\sum_{n=0}^{\infty}\frac{n+1}{(2n+3)(2n+5)}r^{2n}\right]\\
=&\frac{\pi}{6}-\frac{2\pi}{15}=\frac{\pi}{30}.
\end{align*}

(3) Let $f_{12}$ and $f_{14}$ be as in Lemma \ref{Lem1}. Differentiation gives
\begin{align}
r^4f_{23}^{\,\prime}(r)=&2r^2\left\{\left[r-\left(1+r^2\right)\arth r\right]\K+\left(\E-r^{\,\prime2}\K\right)\arth r\right\}\nonumber\\
&+\left[\left(3-r^2\right)\arth r-3r\right]\E+\left(r-r^{\,\prime2}\arth r\right)(\E-\K)\nonumber\\
=&2\left[\left(1+r^2\right)\arth r-r\right]\E+\left[\left(1-5r^2\right)\arth r-r+2r^3\right]\K\nonumber\\
=&2\left[\left(1+r^2\right)\arth r-r\right]\E+r^3f_{12}(r)\K.\label{f22'}
\end{align}
It follows from (\ref{f12}), (\ref{Ineq3}), (\ref{f22'}) and Lemma \ref{Lem1}(10) that
\begin{align*}
r^4f_{23}^{\,\prime}(r)<&2\left[\left(1+r^2\right)\arth r-r\right]\E+2r^3f_{12}(r)\frac{\E}{2-r^2}\\
=&-\frac{2r\E}{2-r^2}\left[3r^{\thinspace\prime2}-\left(3-4r^2-r^4\right)\frac{\arth r}{r}\right]\\
=&-\frac{2r^5\E}{2-r^2}f_{14}(r)<-\frac{52r^5\E}{15(2-r^2)}<0,
\end{align*}
yielding the monotonicity of $f_{23}$.

By Lemma \ref{Lem1}(1), we obtain the limiting values $f_{23}(0^+)=4\pi/3$ and $f_{23}(1^-)=1$.

(4) For $r\in(0,1)$, let
\begin{align*}
G_{13}(r)&=4\left(\E-r^{\thinspace\prime2}\K\right)\frac{\arth r}{r}-3\left(1-\frac{r^{\thinspace\prime2}\arth r}{r}\right)\K \mbox{~and}\\
G_{14}(r)&=r^2\left(1-\frac{r^{\thinspace\prime2}\arth r}{r}\right)\K.
\end{align*}
Then $f_{24}(r)=G_{13}(r)/G_{14}(r)$, $G_{13}(0^+)=G_{14}(0^+)=0$, $G_{14}(1^-)=\infty$, and by (\ref{lim}),
$$G_{13}(1^-)=\lim_{r\to1}\K\left[4\left(\E-r^{\,\prime2}\K\right)\frac{\arth r}{r\K}+3\frac{r^{\,\prime2}\arth r}{r}-3\right]=\infty.$$ Differentiation gives
\begin{align*}
G_{13}'(r)&=\frac{1}{r^2 r^{\thinspace\prime2}}\left[\left(r-r^{\thinspace\prime2}\arth r\right)\E-2r^{\thinspace\prime2}(\arth r-r)\K\right],\\
G_{14}'(r)&=\frac{1}{r^{\thinspace\prime2}}\left[\left(r-r^{\thinspace\prime2}\arth r\right)\E+2r^2r^{\thinspace\prime2}\K\arth r\right], \mbox{~and}\\
\frac{G_{13}'(r)}{G_{14}'(r)}&=\frac{\left(r-r^{\thinspace\prime2}\arth r\right)\E-2r^{\thinspace\prime2}(\arth r-r)\K}{r^2\left[\left(r-r^{\thinspace\prime2}\arth r\right)\E+2r^2r^{\thinspace\prime2}\K\arth r\right]}=\frac{f_{22}(r)}{f_{23}(r)}
\end{align*}
which is strictly increasing on $(0,1)$ by parts (2) and (3). Hence the monotonicity of $f_{24}$ follows from \cite[Theorem 1.25]{AVVb}.

By l'H\^opital's rule, we obtain the limiting values $f_{24}(0^+)=f_{22}(0^+)/f_{23}(0^+)=1/40$ and $f_{24}(1^-)=f_{22}(1^-)/f_{23}(1^-)=1$.

The first lower bound and the upper bound in (\ref{f24}) are clear. The second lower bound in (\ref{f24}) follows from the convexity of $f_{20}$. $\Box$

\begin{corollary}\label{Col1}
The function
$$f_{25}(r)\equiv\frac{[\E(r)-r^{\thinspace\prime2}\K(r)]\arth r}{\left(r-r^{\thinspace\prime2}\arth r\right)\K(r)}$$
is strictly increasing from $(0,1)$ onto $(3/4,1)$. In particular, for $r\in(0,1)$,
\begin{align}\label{Ineqf25}
1-\frac{r-r^{\thinspace\prime2}\arth r}{r^2\arth r}<\frac{\K(r)-\E(r)}{r^2\K(r)}<1-\frac34\frac{r-r^{\thinspace\prime2}\arth r}{r^2\arth r}.
\end{align}
The coefficients 1 and $3/4$ of $(r-r^{\thinspace\prime2}\arth r)/(r^2\arth r)$ in (\ref{Ineqf25}) are both best possible.
\end{corollary}

{\it Proof.} Let $f_{24}$ be as in Theorem \ref{Th5}(4). Then $f_{25}(r)$ can be rewritten as
$$f_{25}(r)=\frac14r^2f_{24}(r)+\frac34,$$
and hence the result for $f_{25}$ follows from Theorem \ref{Th5}(4). The remaining conclusions are clear. $\Box$

\section{\normalsize Proofs of Theorems \ref{Th1}--\ref{Th3}}\label{Sec4}

In this section, we prove Theorems \ref{Th1}--\ref{Th3} stated in Section \ref{Sec1}.

\bigskip
{\it 3.1 Proof of Theorem \ref{Th1}}

\bigskip
 Let $f_2$, $f_7$ and $f_8$ be as in Lemma \ref{Lem1}, $f_{24}$ as in Theorem \ref{Th5}(4), and let
 $$H_1(r)=\frac{f_{24}(r)}{2f_8(r)+1} \mbox{~ and ~} H_2(r)=c-\frac14H_1(r).$$
It follows from Lemma \ref{Lem1}(8) and Theorem \ref{Th5}(4) that $H_1$ is strictly increasing from $(0,1)$ onto $(1/80,1)$. Differentiation gives
\begin{align}\label{1g'}
g'(r)=&2cr\log\frac{\arth r}{r}+\left(\frac34+cr^2\right)\frac{r-r^{\thinspace\prime2}\arth r}{rr^{\thinspace\prime2}\arth r}-\frac{\E-r^{\thinspace\prime2}\K}{rr^{\thinspace\prime2}\K}\nonumber\\
=&cr\left(2\log\frac{\arth r}{r}+\frac{r-r^{\thinspace\prime2}\arth r}{r^{\thinspace\prime2}\arth r}\right)\nonumber\\
&-\frac{4(\E-r^{\thinspace\prime2}\K)\arth r-3(r-r^{\thinspace\prime2}\arth r)\K}{4rr^{\thinspace\prime2}\K\arth r}\nonumber\\
=&cr\left(2\log\frac{\arth r}{r}+\frac{r-r^{\thinspace\prime2}\arth r}{r^{\thinspace\prime2}\arth r}\right)-r\frac{r-r^{\thinspace\prime2}\arth r}{4r^{\thinspace\prime2}\arth r}f_{24}(r)\nonumber\\
=&r\left(2\log\frac{\arth r}{r}+\frac{r-r^{\thinspace\prime2}\arth r}{r^{\thinspace\prime2}\arth r}\right)H_2(r)\nonumber\\
=&r^3\left[f_2(r)+2f_7(r)\right]H_2(r).
\end{align}

By Lemma \ref{Lem1}(8) and Theorem \ref{Th5}(4), $H_2$ is strictly decreasing from $(0,1)$ onto $(c-1/4, c-1/320)$. Hence
it follows from and (\ref{1g'}) that $g$ is strictly increasing (decreasing) on $(0,1)$ if and only if $c\geq1/4$ ($c\leq1/320$, respectively).

Clearly, $g(0^+)=0$. By (\ref{lim}),
\begin{align*}
g(1^-)=&\lim_{r\to1}\log\left(\frac{\pi}{2\K}\left(\frac{\arth r}{r}\right)^{3/4+cr^2}\right)\\
=&\log \frac{\pi}{2}+\lim_{r\to1}\log\left(\frac{\arth r}{r\K}\left(\frac{\arth r}{r}\right)^{cr^2-1/4}\right)\\
=&\log \frac{\pi}{2}+\lim_{r\to1}\left(cr^2-\frac14\right)\log\frac{\arth r}{r}\\
=&\begin{cases}
\log(\pi/2), \mbox{~if} ~c=1/4,\\
\infty, ~~~~~~~~~~~\mbox{if} ~c>1/4,\\
-\infty, ~~~~~~~~\,\mbox{if} ~c<1/4.
\end{cases}
\end{align*}

Next, let $f_{17}$ be as in Lemma \ref{Lem2}(3), and for $r\in(0,1)$, let
\begin{align*}
x&=r^2, ~H_3(r)=g(r)|_{c=1/4},\\
H_4(x)&=H_3\left(\sqrt{x}\right)=\frac{3+x}{4}\log F_1(x)-\log F_0(x),\\
H_5(x)&=(1-x)^{1/4}F_0(x),  ~H_6(x)=(1-x)^{1/3}F_1(x)
\end{align*}
and $H_7(x)=(1-x)F_0(x)$. Then by differentiation and (\ref{F0'F1'}),
\begin{align*}
H_3^{\,\prime}(r)=&2r\left[\frac14\log F_1(x)+\frac{(3+x)F_1'(x)}{4F_1(x)}-\frac{F_0'(x)}{F_0(x)}\right]\\
=&2r\left\{\frac14\log F_1(x)+\frac{F_1'(x)}{F_0(x)F_1(x)}\left[\frac{3+x}{4}F_0(x)-F_1(x)\frac{F_0'(x)}{F_1'(x)}\right]\right\}\\
=&\frac{r}{2}\left\{\log F_1(x)+\frac{G_1(x)}{3(1-x)F_0(x)F_1(x)}\left[4f_{17}(x)-(1-x)F_0(x)\right]\right\}\\
=&\frac{r}{2}\left\{\log F_1(x)+\frac{G_1(x)}{3(1-x)^{5/12}H_5(x)H_6(x)}\left[4f_{17}(x)-H_7(x)\right]\right\},
\end{align*}
which gives
\begin{align}\label{H3'}
\frac{2}{r}H_3^{\,\prime}(r)=H_8(x)\equiv\log F_1(x)+\frac{G_1(x)\left[4f_{17}(x)-H_7(x)\right]}{3(1-x)^{5/12}H_5(x)H_6(x)}.
\end{align}

Clearly, $H_8(0^+)=0$. By \cite[Lemma 2.15(1)]{QVu3}, the functions $H_5$, $H_6$ and $H_7$ are all strictly decreasing from $(0,1)$ onto itself, and by Lemma \ref{Lem2}(3), the function $x\mapsto 4f_{17}(x)-H_7(x)$ is strictly increasing from $(0,1)$ onto $(0, 4(\log4)/\pi)$. Hence $H_8$ is positive and strictly increasing on $(0,1)$, and so is $H_3^{\,\prime}$ by (\ref{H3'}). This shows that $g$ is convex on $(0,1)$ in the case when $c=1/4$.

Let $H_9(r)=r^2\log((\arth r)/r)$. Then
\begin{align}\label{H9'}
H_9^{\,\prime}(r)=r\frac{r-r^{\,\prime2}\arth r}{r^{\,\prime2}\arth r}\left[1+2f_8(r)\right]=r^3\left[f_2(r)+2f_7(r)\right],
\end{align}
which is strictly increasing from $(0,1)$ onto $(0,\infty)$ by Lemma \ref{Lem1}(2) and (7), so that $H_9$ is strictly increasing and convex on $(0,1)$. If $c>1/4$, then $g(r)$ can be written as
\begin{align*}
g(r)=H_3(r)+\left(c-\frac14\right)H_9(r),
\end{align*}
which is a sum of two strictly increasing and convex functions. Hence, if $c\ge1/4$, then $g$ is convex on $(0,1)$.

Finally, if $c\leq1/320$, then by (\ref{1g'}), we have
\begin{align}\label{-g'}
-g'(r)&=r^3\left(\frac{2}{r^2}\log\frac{\arth r}{r}+\frac{r-r^{\,\prime}\arth r}{r^2r^{\,\prime}\arth r}\right)\left[\frac14H_1(r)-c\right]\nonumber\\
&=r^3\cdot\left[f_2(r)+2f_7(r)\right]\cdot\left[\frac14H_1(r)-c\right],
\end{align}
which is a product of three positive and strictly increasing functions on $(0,1)$ by Lemma \ref{Lem1}(2) and (7), and by the monotonicity of $H_1$. Hence $g$ is concave on $(0,1)$ provided that $c\leq1/320$. $\Box$

\bigskip
{\it 3.2 Proof of Theorem \ref{Th2}}

\bigskip
Let $H_1$, $H_3$, $H_8$ and $H_9$ be as in the Proof of Theorem \ref{Th1}.

(1) Clearly, $g_1(r)=H_3(r)/r^2$. Since
\begin{align}\label{H8}
H_3^{\,\prime}(r)\left[\frac{d}{dr}\left(r^2\right)\right]^{-1}=\frac{H_3^{\,\prime}(r)}{2r}=\frac14H_8(x)
\end{align}
by (\ref{H3'}), the monotonicity of $g_1$ follows from that of $H_8$ and \cite[Theorem 1.25]{AVVb}.

Clearly, $g_1(1^-)=H_3(1^-)=\log(\pi/2)$. By l'H\^opital's rule and (\ref{H8}), we obtain
\begin{align*}
g_1(0^+)=\lim_{r\to0}\frac{H_3(r)}{r^2}=\frac14H_8(0^+)=0.
\end{align*}

Observe that $g_2(r)=H_3(r)/H_9(r)$ with $H_3(0^+)=H_9(0^+)=0$. It follows from (\ref{1g'}) and (\ref{H9'}) that
\begin{align}\label{H3'/H9'}
H_3^{\,\prime}(r)/H_9^{\,\prime}(r)=[1-H_1(r)]/4,
\end{align}
which is strictly decreasing from $(0,1)$ onto $(0,79/320)$ by the monotonicity property of $H_1$. Hence the monotonicity of $g_2$ follows from \cite[Theorem 1.25]{AVVb}.

By l'H\^opital's rule and (\ref{H3'/H9'}), $g_2(0^+)=79/320$. Since $H_3(1^-)=\log(\pi/4)$ by part (1), and since $H_9(1^-)=\infty$, we obtain the limiting value $g_2(1^-)=0$.

(2) For $r\in(0,1)$, let
$$H_{10}(r)=-g(r)|_{c=1/320}=\log\frac{2\K(r)}{\pi}-\left(\frac34+\frac{1}{320}r^2\right)\log\frac{\arth r}{r}.$$
Then $H_9(0^+)=H_{10}(0^+)=0$, $g_3(r)=H_{10}(r)/H_9(r)$, and it follows from (\ref{H9'}) and (\ref{-g'}) that
\begin{align}\label{H10'H9'}
\frac{H_{10}^{\,\prime}(r)}{H_9^{\,\prime}(r)}=\frac{1}{4}\left[H_1(r)-\frac{1}{80}\right],
\end{align}
which is strictly increasing on $(0,1)$ by the monotonicity property of $H_1$. Hence the monotonicity of $g_3$ follows from \cite[Theorem 1.25]{AVVb}.

By l'H\^opital's rule, (\ref{H10'H9'}), and Corollary \ref{Col1}, we obtain the liming values $g_3(0^+)=0$ and
\begin{align*}
g_3(1^-)=&\lim_{r\to1}\frac{\log(2\K/\pi)}{\log((\arth r)/r)}-\frac34-\frac{1}{320}\\
=&\lim_{r\to1}\frac{[\E(r)-r^{\thinspace\prime2}\K(r)]\arth r}{\left(r-r^{\thinspace\prime2}\arth r\right)\K(r)}-\frac{241}{320}=\frac{79}{320}.
\end{align*}

Next, let $f_7$ be as in Lemma \ref{Lem1}. Then
$$g_4(r)=\left(\frac{1}{r^2}\log\frac{\arth r}{r}\right)g_3(r)=f_7(r)g_3(r),$$
and hence the result for $g_4$ follows from Lemma \ref{Lem1}(7) and the property of $g_3$.

(3) It is clear that
$$f(r)=\frac14-g_2(r) \mbox{~ and ~} G(r)=r^2f(r)+\frac34.$$
Hence the results for $f$ and $G$ follow from part (1).

The assertion on (\ref{Ineq1}) is obvious. The double inequality (\ref{Ineq2}) follows from  the property of $g_1$. $\Box$

\bigskip
{\it 3.3 Proof of Theorem \ref{Th3}}

\bigskip
(1) The monotonicity and limiting values of $h_1$ were obtained in \cite[Theorem 3.10]{AVV2}. By differentiation,
\begin{align}\label{h1'}
h_1'(r)=\frac{\E(r)\arth r-r\K(r)}{r^{\thinspace\prime2}\arth^2 r}=-H_{11}(r)H_{12}(r),
\end{align}
where
$$H_{11}(r)=r\left(\frac{r}{r^{\thinspace\prime}\arth r}\right)^2 \mbox{~ and ~} H_{12}(r)=\frac{r\K(r)-\E(r)\arth r}{r^3}.$$
By \cite[Lemma 3(1)]{QVV1}, $H_{11}$ is strictly increasing from $(0,1)$ onto $(0,\infty)$. Let $H_{13}(r)=r\K(r)-\E(r)\arth r$ and $H_{14}(r)$. Then $H_{13}(0)=H_{14}(0)=0$, $H_{12}(r)=H_{13}(r)/H_{14}(r)$, and
$$\frac{H_{13}'(r)}{H_{14}'(r)}=\frac{\K(r)-\E(r)}{3r^2}\cdot\frac{\arth r}{r}$$
which is strictly increasing from $(0,1)$ onto $(\pi/12,\infty)$ by \cite[Lemma 5.2 (3)]{AQVV}, so that $H_{12}$ is positive and strictly increasing on $(0,1)$ by \cite[Theorem 1.25]{AVVb}. Consequently, it follows from (\ref{h1'}) that $h_1'$ is strictly decreasing on $(0,1)$. This yields the concavity of $h_1$.

(2) Let $c_n$ be as in Lemma \ref{Lem4}. Then by Lemma \ref{Lem4}(1),
\begin{align}\label{tildeak-ck}
\tilde{a}_{k+1}=\frac{1}{2k+3}\left[\frac34-(2k+3)a_{k+1}\right]=\frac{1}{2k+3}\left(\frac34-c_k\right)>0
\end{align}
for $k\in\IN_0$. It follows from (\ref{F})--(\ref{arth}) that
\begin{align}\label{h2}
h_2(r)
=&\frac{1}{r^{2(n+1)}}\left(\frac34\sum_{n=0}^{\infty}\frac{1}{2n+1}r^{2n}-\sum_{n=0}^{\infty}a_nr^{2n}-\sum_{k=0}^n\tilde{a}_k r^{2k}\right)\nonumber\\
=&\frac{1}{r^{2(n+1)}}\left(\sum_{n=0}^{\infty}\tilde{a}_n r^{2n}-\sum_{k=0}^n\tilde{a}_k r^{2k}\right)
=\sum_{k=0}^{\infty}\tilde{a}_{k+n+1}r^{2k}.
\end{align}
Hence by (\ref{tildeak-ck}), all the coefficients of the Maclaurin series of $h_2$ are positive, so that $h_2$ is strictly absolutely monotone on $(0,1)$.

By (\ref{tildeak-ck}) and (\ref{h2}), $h_2(0^+)=\tilde{a}_{n+1}>0$, and by (\ref{Form2}), we obtain
\begin{align*}
h_2(1^-)
=&\lim_{r\to1}\left(\frac{3}{8r}\log\frac{1+r}{1-r}-\frac{2}{\pi}\log\frac{4}{r^{\,\prime}}\right)-\sum_{k=0}^{n}\tilde{a}_k\\
=&\lim_{r\to1}\left(\frac{3}{8r}-\frac{1}{\pi}\right)\log\frac{1}{1-r}-\sum_{k=0}^{n}\tilde{a}_k-\frac{8-\pi}{8\pi}\log8=\infty.
\end{align*}

(3) It follows from (\ref{arth}) and (\ref{h2})
\begin{align}\label{h3}
h_3(r)=\left(\sum_{k=0}^{\infty}\tilde{a}_{k+n+1}r^{2k}\right)\left(\sum_{k=0}^{\infty}\frac{1}{2k+1}r^{2k}\right)^{-1}.
\end{align}
The ratio of the coefficients of the two series expansions in (\ref{h3}) is equal to
\begin{align}\label{ratioCoeffh3}
\frac{2k+1}{2(k+n)+3}\left(\frac{3}{4}-c_{k+n}\right)
=\left[1-2\frac{n+1}{2(k+n)+3}\right]\left(\frac{3}{4}-c_{k+n}\right).
\end{align}
By Lemma \ref{Lem4}(1), the left side of (\ref{ratioCoeffh3}) is a product of two factors which are both positive and strictly increasing in $k\in\IN_0$. Hence the monotonicity of $h_3$ follows from \cite[Lemma 2.1]{PV}.

Clearly, $h_3(0^+)=h_2(0^+)=\tilde{a}_{n+1}$. By part (1), we obtain the limiting value
\begin{align*}
h_3(1^-)=\lim_{r\to1}\left[\frac34-\frac{2}{\pi}h_1(r)-\frac{r}{\arth r}\sum_{k=0}^{n}\tilde{a}_k\right]=\frac34-\frac{2}{\pi}.
\end{align*}

(4) Let $d_n$ be as in Lemma \ref{Lem4}. Then by (\ref{F})--(\ref{arth}), we have
\begin{align*}
&\frac{\pi}{2}\left(1-\frac{1}{12}r^2-\frac{91}{2880}r^4\right)\frac{\arth r}{r}-\K(r)\\
&=\frac{\pi}{2}\left(\sum_{n=1}^{\infty}\frac{r^{2n}}{2n+1}-\frac{1}{12}r^2\sum_{n=0}^{\infty}\frac{r^{2n}}{2n+1}
-\frac{91}{2880}r^4\sum_{n=0}^{\infty}\frac{r^{2n}}{2n+1}-\sum_{n=1}^{\infty}a_nr^{2n}\right)\\
&=\frac{\pi}{2}\left(\sum_{n=2}^{\infty}\frac{r^{2n}}{2n+1}-\frac{1}{12}r^2\sum_{n=1}^{\infty}\frac{r^{2n}}{2n+1}
-\frac{91}{2880}r^4\sum_{n=0}^{\infty}\frac{r^{2n}}{2n+1}-\sum_{n=2}^{\infty}a_nr^{2n}\right)\\
&=\frac{\pi}{2}r^4\sum_{n=0}^{\infty}\left[\frac{1}{2n+5}-\frac{1}{12(2n+3)}-\frac{91}{2880(2n+1)}-a_{n+2}\right]r^{2n}\\
&=\frac{\pi}{5760}r^4\sum_{n=0}^{\infty}\left[\frac{10196n^2+18704n+6075}{(2n+5)(2n+3)(2n+1)}-2880a_{n+2}\right]r^{2n}\\
&=\frac{\pi}{5760}r^6\sum_{n=0}^{\infty}\left[\frac{10196n^2+39096n+34975}{(2n+7)(2n+5)(2n+3)}-2880a_{n+3}\right]r^{2n}\\
&=\frac{\pi}{5760}r^6\sum_{n=0}^{\infty}\frac{d_n}{2n+1}r^{2n},
\end{align*}
and hence
\begin{align}\label{h4}
h_4(r)=&\frac{1}{(\arth r)/r}\cdot\frac{1}{r^6}\left[\frac{\pi}{2}\left(1-\frac{1}{12}r^2-\frac{91}{2880}r^4\right)\frac{\arth r}{r}-\K(r)\right]\nonumber\\
=&\frac{\pi}{5760}\left(\sum_{n=0}^{\infty}\frac{d_n}{2n+1}r^{2n}\right)
\cdot\left(\sum_{n=0}^{\infty}\frac{1}{2n+1}r^{2n}\right)^{-1}.
\end{align}
Clearly, the ratio of the coefficients of the two series expansions in (\ref{h4}) is  $d_n$.
Hence the monotonicity of $h_4$ follows from (\ref{h4}), Lemma \ref{Lem4}(3) and \cite[Lemma 2.1]{PV}.

By (\ref{h4}), $h_4(0^+)=\pi d_0/5760=871\pi/96768$, and by part (1),
$$h_4(1^-)=\frac{\pi}{2}\left(1-\frac{1}{12}-\frac{91}{2880}\right)-h_1(1^-)=\frac{2549\pi}{5760}-1.$$

(5) The double inequality (\ref{K-arth1}) ((\ref{K-arth2})) follows from part (3) (part (4), respectively). Taking $n=1$ in (\ref{K-arth1}), we obtain the double inequality (\ref{K-arth3}).

By parts (3)--(4), the remaining conclusion is clear. $\Box$

\section{\normalsize Concluding Remark}\label{Sec5}

(i) The double inequality (\ref{Ineq1}) can also be obtained from Theorem \ref{Th1}, or Theorem \ref{Th2}(1), or Theorem \ref{Th2}(2).

(ii) The upper bound of $\K(r)$ given in (\ref{Ineq1}) coincides with that given in (\ref{KSIneq}) obtained by Kup\'an and Sz\'asz. However, our result is implied by a stronger conclusion, and our proof is more direct and simpler than theirs. In addition, it was not indicated in \cite{KS} that the constant $\beta=1/4$ in (\ref{Ineq1}) is best possible.

On the other hand, since the function
\begin{align*}
h_5(r)\equiv\left\{\left(\frac{\pi}{2}\right)^{1-r^2}\left(\frac{\arth r}{r}\right)^{3/4+r^2/4}\left[\frac{\pi}{2}\left(\frac{\arth r}{r}\right)^{3/4+r^2/320}\right]^{-1}\right\}^{1/r^2}=\frac{2}{\pi}\left(\frac{\arth r}{r}\right)^{79/320}
\end{align*}
is strictly increasing from $(0,1)$ onto $(2/\pi,\infty)$, there exists a number $r_0\in(0,1)$ satisfying the condition $(\arth r_0)/r_0=(\pi/2)^{320/79}$ such that the lower bound given in (\ref{Ineq1}) ( (\ref{Ineq2}) ) is better than that given in (\ref{Ineq2}) ( (\ref{Ineq1}) ) for $r\in(0,r_0)$ ( $r\in(r_0,1)$, respectively).

(iii) Let $h_4$ be as in Theorem \ref{Th3}. Then by Theorem \ref{Th3}(4), the function
$$h_6(r)\equiv \frac{r\K(r)}{\arth r}+\frac{\pi}{24}r^2+\frac{91\pi}{5760}r^4=\frac{\pi}{2}-r^6h_4(r)$$
is strictly decreasing from $(0,1)$ onto $(1+331\pi/5760,\pi/2)$. Therefore, Theorem \ref{Th3}(4) remarkably improves the monotonically decreasing property of $h_1(r)=r\K(r)/\arth r$ obtained in \cite[Theorem 3.10]{AVV2}.

(iv) It follows from Theorem \ref{Th3}(1) that for all $r\in(0,1)$,
\begin{align}\label{Bound1OfK}
\left[1-\left(1-\frac{2}{\pi}\right)r\right]\frac{\arth r}{r}<\frac{2}{\pi}\K(r)<\frac{\arth r}{r}.
\end{align}
Clearly, each of the upper bounds of $2\K(r)/\pi$ given in (\ref{Ineq1})--(\ref{Ineq2}) and in (\ref{K-arth2})--(\ref{K-arth3}) is better than the upper bound given in (\ref{Bound1OfK}).

Let $P_3$, $\alpha$ and $\delta$ be as in Theorem \ref{Th3}(5), and for $r\in(0,1)$, let
\begin{align*}
h_7(r)&=\frac{1}{r}\left\{\left[P_3(r)-\alpha r^6\right]-\left[1-\left(1-\frac{2}{\pi}\right)r\right]\right\} \mbox{~ and}\\
h_8(r)&=\frac{1}{r}\left\{\left[\frac14+\frac34\left(1-\delta r^4\right)\frac{\arth r}{r}\right]-\left[1-\left(1-\frac{2}{\pi}\right)r\right]\frac{\arth r}{r}\right\}.
\end{align*}
Then
\begin{align}
h_7(r)&=1-\frac{2}{\pi}-\frac{1}{12}r-\frac{91}{2880}r^3-\alpha r^5 \mbox{~ and} \label{h7}\\
h_8(r)&=\frac{1}{r}\left\{\frac14+\left[\left(1-\frac{2}{\pi}\right)r-\left(\frac34-\frac{2}{\pi}\right)r^4-\frac14\right]\frac{\arth r}{r}\right\}\nonumber\\
&=\frac{1}{4r}\left[1-\frac{(1-r)\arth r}{r}\right]+\left(\frac34-\frac{2}{\pi}\right)\left(1-r^3\right)\frac{\arth r}{r},\label{h8}
\end{align}
from which we can easily see that $h_7$ is strictly decreasing from $(0,1)$ onto $(0, 1-2/\pi)$, and $h_8(r)>0$ for $r\in(0,1)$ since the function
$$r\mapsto \frac{(1-r)\arth r}{r}=\frac{1}{1+r}\frac{r^{\thinspace\prime2}\arth r}{r}$$
is strictly decreasing from $(0,1)$ onto itself by \cite[Lemma 3(1)]{QVV1}. Consequently, each of the lower bounds of $2\K(r)/\pi$ given in (\ref{K-arth2})--(\ref{K-arth3}) is better than the lower bound of this function given in (\ref{Bound1OfK}).

(v) The pair of lower bounds, and the pair of upper bounds of $\K(r)$ given in (\ref{K-arth2})--(\ref{K-arth3}) are both not directly comparable on $(0,1)$. In order to explain this, we let
\begin{align*}
h_9(r)=&\left[1-\frac{1}{12}r^2-\frac{91}{2880}r^4-\left(\frac{2549}{2880}-\frac{2}{\pi}\right)r^6\right]\frac{\arth r}{r}-\left\{\frac14+\frac34\left[1-\left(1-\frac{8}{3\pi}\right)r^4\right]\frac{\arth r}{r}\right\}\\
=&\left[\frac14-\frac{1}{12}r^2+\left(\frac{2069}{2880}-\frac{2}{\pi}\right)r^4-\left(\frac{2549}{2880}-\frac{2}{\pi}\right)r^6\right]\frac{\arth r}{r}-\frac14,\\
h_{10}(r)=&\left(1-\frac{1}{12}r^2-\frac{91}{2880}r^4-\frac{871}{48384}r^6\right)\frac{\arth r}{r}-\left[\frac14+\left(\frac34-\frac{3}{320}r^4\right)\frac{\arth r}{r}\right]\\
=&\frac14\left(\frac{\arth r}{r}-1\right)-r^2\left(\frac{1}{12}+\frac{1}{45}r^2+\frac{871}{48384}r^4\right)\frac{\arth r}{r}.
\end{align*}
Applying l'H\^opital's, one can easily obtain the limiting value
\begin{align*}
\lim_{r\to1}\frac{1}{r^{\thinspace\prime2}}\left[\frac14-\frac{1}{12}r^2+\left(\frac{2069}{2880}-\frac{2}{\pi}\right)r^4
-\left(\frac{2549}{2880}-\frac{2}{\pi}\right)r^6\right]=\frac{1}{12},
\end{align*}
and hence $h_9(1^-)=-1/4$. On the other hand, by (\ref{arth}), we obtain
\begin{align*}
\frac{3}{r^4}h_9(r)=\sum_{n=0}^{\infty}\left[\frac{n+1}{(2n+3)(2n+5)}+\left(\frac{2069}{960}-\frac{6}{\pi}\right)\frac{1}{2n+1}\right]r^{2n}
-\left(\frac{2549}{960}-\frac{6}{\pi}\right)\sum_{n=0}^{\infty}\frac{r^{2(n+1)}}{2n+1},
\end{align*}
so that
\begin{align*}
\lim_{r\to0}\frac{3}{r^4}h_9(r)=\frac{1}{15}+\frac{2069}{960}-\frac{6}{\pi}=\frac{2133}{960}-\frac{6}{\pi}=0.3120156\cdots>0.
\end{align*}
Consequently, there exist two numbers $r_1, r_2\in(0,1)$ such that $h_9(r)$ is positive (negative) for $r\in(0,r_1)$ ($r\in(r_2,1)$, respectively). This shows that the lower bounds of $\K(r)$ given in (\ref{K-arth2})--(\ref{K-arth3}) are not directly comparable on $(0,1)$.

Similarly, we have
\begin{align*}
\frac{241920}{r^6}h_{10}(r)=\sum_{n=0}^{\infty}\frac{244712n^3+749388n^2+408406n-161595}{(2n+1)(2n+3)(2n+5)(2n+7)}r^{2n},
\end{align*}
by which we obtain $\lim_{r\to0}[241920 h_{10}(r)/r^6]=-161595/105=-1539$. On the other hand, we have
\begin{align*}
\frac{h_{10}(r)}{r\arth r}=\frac{\arth r-r}{4r^2\arth r}-\left(\frac{1}{12}+\frac{1}{45}r^2+\frac{871}{48384}r^4\right),
\end{align*}
from which we see that
\begin{align*}
\lim_{r\to1}\frac{h_{10}(r)}{r\arth r}=\frac14-\frac{1}{12}-\frac{1}{45}-\frac{871}{48384}=\frac{30589}{241920}>0.
\end{align*}
Therefore, there exist two numbers $r_3, r_4\in(0,1)$ such that $h_{10}(r)$ is negative (positive) for $r\in(0,r_3)$ ($r\in(r_4,1)$, respectively),  showing that the upper bounds of $\K(r)$ given in (\ref{K-arth2})--(\ref{K-arth3}) are also not directly comparable on $(0,1)$.

(vi) Let $f$, $G$, $g_1$ and $g_4$ be as in Theorem \ref{Th2}.
By computation, we have
\begin{align}
f(r)=&\frac{1}{320}+\frac{517}{201600}r^2+\frac{767341}{387072000}r^4+\frac{4277471797}{2682408960000}r^6\nonumber\\
&+\frac{1851483120061}{1394852659200000}r^8+\frac{2989339649544551}{2636271525888000000}r^{10}+\cdots,\label{Serf}\\
G(r)=&\frac34+\frac{1}{320}r^2+\frac{517}{201600}r^4+\frac{767341}{387072000}r^6+\frac{4277471797}{2682408960000}r^8\nonumber\\
&+\frac{1851483120061}{1394852659200000}r^{10}+\frac{2989339649544551}{2636271525888000000}r^{12}+\cdots,\label{SerG}\\
h_{11}(r)\equiv &\frac{g_1(r)}{r^2}=\frac{1}{r^4}\left[\frac{3+r^2}{4}\log\frac{\arth r}{r}-\log\frac{2\K(r)}{\pi}\right]\nonumber\\
=&\frac{79}{960}+\frac{421}{12096}r^2+\frac{690961}{33177600}r^4+\frac{18414493}{1277337600}r^6\nonumber\\
&+\frac{164673431213}{15216574464000}r^8+\cdots\label{Ser1}
\end{align}
and
\begin{align}
h_{12}(r)\equiv &\frac{g_4(r)}{r^2}=\frac{1}{r^6}\left[\log\frac{2\K(r)}{\pi}-\left(\frac34+\frac{1}{320}r^2\right)\log\frac{\arth r}{r}\right]\nonumber\\
=&\frac{517}{604800}+\frac{239497}{232243200}r^2+\frac{741527}{709632000}r^4\nonumber\\
&+\frac{168874886801}{167382319104000}r^6
+\frac{2405137262477}{2510734786560000}r^8+\cdots\label{Ser2}
\end{align}
Our computation seems to show that the following conjectures are true.
\begin{conjecture}\label{Conj}
(1) The function $f$ in Theorem \ref{Th2}(3) is convex on $(0,1)$.

(2) The functions $h_{11}$ and $h_{12}$ are both strictly increasing and convex on $(0,1)$, with ranges $(79/960, \log(\pi/2))$ and $(517/604800, \infty)$, respectively.

(3) Let $f_7$ be as in Lemma \ref{Lem1}(7), and let $h_{13}(r)=r^{-2}\log(2\K(r)/\pi)$ for $r\in(0,1)$. Then the coefficients of the Maclaurin series expansions of the functions $f$, $G$, $f_7$, $h_{11}$, $h_{12}$ and $h_{13}$ are all positive, and in particular, these functions are all absolutely monotone on $(0,1)$.
\end{conjecture}

If these conjectures are true, then the inequalities in Theorem \ref{Th2}(3) can be improved. For example, the following double inequalities
\begin{align}\label{IneqConj1}
\left(\frac{\arth r}{r}\right)^{\frac34+\frac{1}{320}r^2+\frac{517}{201600}r^4}<\frac{2\K(r)}{\pi}<\left(\frac{\arth r}{r}\right)^{\frac34+\frac{1}{320}r^2+\frac{79}{320}r^4}
\end{align}
and
\begin{align}\label{IneqConj2}
\tau^{r^4}\left(\frac{\arth r}{r}\right)^{\frac34+\frac14r^2}<\frac{2\K(r)}{\pi}<\mu^{r^4}\left(\frac{\arth r}{r}\right)^{\frac34+\frac14r^2}
\end{align}
will be valid for all $r\in(0,1)$, where $\tau=2/\pi=0.6366197\cdots$ and $\mu=e^{-79/960}=0.9210032\cdots$.

(vii) We also would like to put forward the following
\begin{problem}\label{P}
What are the analogues of Theorems \ref{Th1}--\ref{Th3} for the generalized complete elliptic integrals (or even for the zero-balanced hypergeometric functions)?
\end{problem}

\bigskip
\noindent {\bf Acknowledgments}

(1) This work is supported by the National Natural Science Foundation of China (Grant No.11771400).

(2) The two universities at which the corresponding author works are listed in no particular order.

\bigskip
\noindent {\bf Competing Interests}

The authors declare that they have no competing interests.


\bigskip
\bigskip
\noindent Authors' addresses:

\bigskip
\noindent Song-Liang Qiu:\\
Department of Mathematics, Lishui University, Lishui 323000, Zhejiang, China;\\
Department of Mathematics, Zhejiang Sci-Tech University, Hangzhou 310018, Zhejiang, China\\
E-mail address: sl$\_$qiu@zstu.edu.cn

\bigskip
\noindent Qi Bao:\\
School of Mathematical Sciences, East China Normal University, Shanghai 200241, China\\
E-mail address: 52205500010@stu.ecnu.edu.cn

\bigskip
\noindent Xiao-Yan Ma:\\
Department of Mathematics, Zhejiang Sci-Tech University, Hangzhou 310018, Zhejiang, China\\
E-mail address: mxy@zstu.edu.cn

\bigskip
\noindent Hong-Biao Jiang:\\
Department of Mathematics, Lishui University, Lishui 323000, Zhejiang, China\\
E-mail address: lsxyhbj@126.com


\begin{thebibliography}{}

%
%
\bibitem{Ask2} R. Askey, in: P. Duren (Ed.), Handbooks of Special Functions, A century of Mathematics in America, Part III, Amer. Math. Soc. (1989)

\bibitem{AAR}
G. Andrews, R. Askey and R. Roy, Special Functions, Encyclopedia of Mathematics and Its Applications, vol. 71, Cambridge Univ Press (1999)

\bibitem{AQ}
H. Alzer and S.L. Qiu, Monotonicity theorems and inequalities for the complete elliptic integrals, J. Comput. Appl. Math. 172(2), 289--312 (2004)

\bibitem{AS}
M. Abramowitz and I.A. Stegun, Handbook of Mathematical Functions With Formulas, Graphs and Mathematical Tables. New York: Dover, 1965. Tenth Printing, with corrections (1972)

\bibitem{AQVa}
G.D. Anderson, S.L. Qiu, and M.K. Vamanamurthy, Elliptic integrals inequalities, with applications, Constr. Approx. 14, 195--207 (1998)

\bibitem{AQVV}
G.D. Anderson, S.L. Qiu, M.K. Vamanamurthy, and M. Vuorinen, Generalized elliptic integrals and modular equations, Pacific J. Math. 192(1), 1--37 (2000)

\bibitem{AVV1}
G.D. Anderson, M.K. Vamanamurthy, and M. Vuorinen, Functional inequalities for complete elliptic ingrals and their ratios, SIAM J. Math. Anal. 21(2), 536--549 (1990)

\bibitem{AVV2}
G.D. Anderson, M.K. Vamanamurthy, and M. Vuorinen, Functional inequalities for hypergeometric functions and complete elliptic integrals, SIAM J. Math. Anal. 21(2), 512--524 (1992)

\bibitem{AVVb}
G.D. Anderson, M.K. Vamanamurthy, M. Vourinen, Conformal Invarinants, Inequalities, and Quasiconformal Maps, John Wiley and Sons, New York (1997)

\bibitem{B}
A. Baricz, Tur\'an type inequalities for generalized complete elliptic integrals, Math. Z. 256(4), 895--911 (2007)

\bibitem{Be2}
B.C. Berndt, Ramanujan's Notebooks, Part II, Springer-Verlag, New York (1989)

\bibitem{BB}
J.M. Borwein and P.B. Borwein, Pi and the AGM, John Wiley \& Sons, New York (1987)

\bibitem{C}
B.C. Carlson, Special Functions of Applied Mathematics, Academic Press, New York (1977)

\bibitem{CWQ1}
Y.-M. Chu, M.-K. Wang and Y.-F. Qiu, On Alzer and Qiu's conjecture for complete elliptic integral and inverse hyperbolic tangent function, Abstract and Applied Analysis 2011, Article ID 697547, 7 pages (https://doi.org/10.1155/2011/697547)

\bibitem{CWQ}
Y.M. Chu, M.K. Wang, S.L. Qiu, and Y.P. Jiang, Bounds for complete
elliptic integrals of the second kind with applications, Comput. Math. Appl. 63(7), 117--1184 (2012)

\bibitem{HQM}
T.R. Huang, S.L. Qiu and X.Y. Ma, Monotonicity properties and inequalities for the generalized  elliptic integral of the first kind, J. Math. Anal. Appl. 469(1), 95--116 (2019)

\bibitem{KS}
P\'al A. Kup\'an and R\'obert Sz\'asz, About bounds for the elliptic integral of the first kind, Rev. Anal. Num\'er. Th\'eor. Approx. 41(2), 149--156 (2012)

\bibitem{OLB}
F.W.J. Olver, D.W. Lozier, R.F. Boisvert, et al., NIST Handbook of Mathematical Functions, Cambridge University Press, Cambridge (2010)

\bibitem{PV}
S. Ponnusany, M. Vourinen, Asymptotic expansions and inequalities for hypergeometric functions, Mathematika 44(2), 278--301 (1997)

\bibitem{QMC2}
S.L. Qiu, X.Y. Ma and Y.M. Chu, Sharp Landen transformation inequalities for hypergeometric functions, with applications, J. Math. Anal. Appl. 474(2), 1306--1337 (2019)

\bibitem{QVa}
 S.L. Qiu and M.K. Vamanamurthy, Sharp estimates for complete elliptic integrals, SIAM
J. Math. Anal. 27(3), 823--834 (1996)

\bibitem{QVV1}
S.L. Qiu, M.K. Vamanamurthy and M. Vuorinen, Some inequalities for the Hersch-Pfluger distortion function.
J. of Inequal. \& Appl. 4, 115--139 (1999)

\bibitem{QVu2}
S.L. Qiu and M. Vuorinen, Handbook of Complex Analysis: Special Function in
Geometric Function Theory, pp.621--659. Elsevier B.V., Amsterdam (2005)

\bibitem{QVu3}
S.L. Qiu and M. Vuorinen, Duplication inequalities for the ratios of hypergeometric functions, Forum Math. 12(1), 109--133 (2000)

\bibitem{WZC}
G.D. Wang, X.H. Zhang, Y.M. Chu and S.L. Qiu, Complete elliptic integrals and the Hersch-Pfluger distortion function (in Chinese), Acta Mathematica Scientia 4, 731--734 (2008)

\end{thebibliography}
\end{document}